\documentclass[12pt]{article}
\usepackage{amssymb,latexsym, amsmath}
\usepackage{color}

\def\rr{{\mathbb R}}
\def\rn{{{\rr}^n}}

\def\zz{{\mathbb Z}}
\def\ls{\lesssim}
\def\N{\mathbb{N}}

\def\R{\mathbb{R}}

\newtheorem{theo}{\textsc{Theorem}}
\newtheorem{lem}[theo]{\textsc{Lemma}}
\newtheorem{propo}[theo]{\textsc{Proposition}}

\title{ Sharp endpoint estimates for some operators associated with the Laplacian with drift in Euclidean space}
\author{Hong-Quan Li and Peter Sj\"ogren}
\date{}
\begin{document}
\renewcommand{\theequation}{\thesection.\arabic{equation}}
\setcounter{equation}{0} \maketitle

\vspace{-1.0cm}

\bigskip

\mbox{}\\
Hong-Quan Li\\
School of Mathematical Sciences,
Fudan University,
220 Handan Road,
Shanghai 200433,
People's Republic of China\\
E-Mail: hongquan\_li@fudan.edu.cn \quad or \quad hong\_quanli@yahoo.fr \\
\mbox{}\\
Peter Sj\"ogren (corresponding author)\\
Mathematical Sciences,
University of Gothenburg, \hskip3pt
 and \hskip3pt \\
Mathematical Sciences,
Chalmers,  \\
 SE-412 96 G\"oteborg,
Sweden\\
E-Mail: peters@chalmers.se \\
Tel.: +46 31 7723511, Fax  +46 31 161973\\

\bigskip

\noindent {\bf Abstract.} Let $v \ne 0$ be a vector in $\R^n$.
 Consider the Laplacian on $\R^n$ with drift
$\Delta_{v} = \Delta + 2v\cdot \nabla$ and the measure $d\mu(x) = e^{2 \langle v, x \rangle}
dx$, with respect to which $\Delta_{v}$ is self-adjoint.
This measure has exponential growth with respect to the Euclidean distance.
We study weak type $(1, 1)$  and other sharp endpoint estimates for the Riesz
transforms of any order, and also for the vertical and horizontal
Littlewood-Paley-Stein
functions  associated with the heat and the Poisson semigroups.

\bigskip


\bigskip


\newpage

\medskip

\renewcommand{\theequation}{\thesection.\arabic{equation}}
\section{Introduction}
\setcounter{equation}{0}
\medskip

Consider the weighted manifold $\R^{(n, v)}$,  
defined as $\R^n$ with the Euclidean distance and the measure $ d\mu(x) =e^{2\langle v,x \rangle} dx$. Here
  $v = (v_1, \cdots, v_n) \in \R^n \setminus \{
0 \}$ is fixed and $\langle ., . \rangle$ denotes  the inner product on $\R^n$.
 Notice that large balls in this space have exponential volume growth.
With $\R^{(n, v)}$, we associate the Laplacian with drift  
\begin{align*}
\Delta_{v} = \sum_{i = 1}^n \Big( \frac{\partial^2}{\partial x_i^2} + 2 v_i \frac{\partial}{\partial x_i} \Big).
\end{align*}

We denote by $\nabla$ the usual gradient operator in $\R^n$.
Notice that the Green formulas holds with respect to the measure $\mu$, that is,
\begin{align*}
\int_{\R^{n}} f \Delta_{v}  w \, d\mu = - \int_{\R^{n}} \langle \nabla f, \nabla w \rangle \, d\mu = \int_{\R^{n}} w \Delta_{v}  f \, d\mu,
\end{align*}
provided  $f$ and $w$ are smooth and $f$ or $w$ has compact support.
Thus $\Delta_{v}$ has a selfadjoint extension in  $L^2(\R^{n}, d\mu)$.

  We shall consider Riesz transforms and Littlewood-Paley-Stein functions of any
order in  $\R^{(n, v)}$. These operators are defined and studied in many general
settings, such as Lie groups, symmetric spaces and other Riemannian manifolds.
Their $L^p$ boundedness properties for $1<p<\infty$  have been well studied and
are known in several cases.
 We refer the reader to \cite{A92},  \cite{ACDH04}-\cite{ERS99}, \cite{Li99}, \cite{L85}-\cite{LM04} and references therein.
In particular, the classical results are not always valid on manifolds.
Interesting counterexamples can be found in \cite{Li99}, \cite{CD99}, \cite{CCH06} and \cite{CCFR16} for Riesz transforms of order one, and in \cite{GQS92}, \cite{A92} and \cite{ERS99} for those of order two.
  The setting of the Ornstein-Uhlenbeck semigroup in Euclidean space is considered in \cite{FGS94} and  \cite{PS00}; see also   \cite{AS12}. There the measure is Gaussian, thus finite but not doubling.

The weak type $(1,1)$ property of these operators is more difficult and less known in spaces of
exponential volume growth. The main reason is that the existing methods to
treat singular integrals are insufficient. Theorem 1.2 of
 \cite{CD99} is a weaker estimate.
 But some results have been established; see \cite{GQS92} and \cite{SV08} for the affine group, which is a typical case without spectral gap. Some other groups and spaces are treated in \cite{An92}, \cite{ADY96} and \cite{HS03}. 
Further,  \cite{LSW14}
and  \cite{LS15} deal with the Laplacian with drift.

Returning to our setting   $\R^{(n, v)}$, we mention that
  Lohou\'e  and Mustapha \cite{LM04} proved that the
 Riesz transforms $\nabla^k (-\Delta_{v} )^{-{k}/{2}}$ of any order $k$
are bounded on $L^p,\; 1<p<\infty$. Their setting and results are
 actually more general.  In \cite{LSW14},
the authors and Y.-R. Wu showed that the first-order Riesz transform
 $\nabla (-\Delta_{v} )^{-{1}/{2}}$ is of weak type $(1, 1)$ in  $\R^{(n, v)}$.
 Here and in the sequel, $L^p$ and weak  $L^p$   estimates in  $\R^{(n, v)}$
always refer to the measure $\mu$.

For   $\R^{(n, v)}$ we observe that ${\partial}/{\partial x_i}$  commutes
with $\Delta_v$ and thus with any negative power of $-\Delta_v$,
so that the factors  $\nabla^k$  and $(-\Delta_{v} )^{-{k}/{2}}$ can be written
in any order.  We  will study  the  weak type $(1, 1)$ property of
  $\nabla^k (-\Delta_{v} )^{-{k}/{2}}$.
 But instead of  $\nabla^k$, we will
 use a  general homogeneous differential operator of order $k\ge 1$,
\begin{equation}\label{D}
D = \sum_{|\alpha|=k} a_\alpha \partial^\alpha
\end{equation}
with constant coefficients, not all 0.
Our Riesz operator will thus be
\begin{equation*}
R_D = D(-\Delta_{v})^{-\frac k2}.
\end{equation*}
 Letting $\partial_v$ denote differentiation along
the vector  $v$, we can write $D$ as a sum
\begin{equation*}
D = \sum_{i=0}^k \partial_v^i D_{k-i}',
\end{equation*}
where $ D_{k-i}'$ is a constant coefficient   operator of order $k-i$ involving
only differentiation in directions orthogonal to   $v$.
 The maximal order of
differentiation along  $v$ is then
\begin{equation*}
q = \max\, \{i: D_{k-i}'\ne 0\} \in \{0,\dots, k\},
\end{equation*}
and this quantity  turns out to be significant.

Our result about $R_D$ is the following.

\begin{theo} \label{one'}
With $D$ and $q$ as just described, the Riesz transform $R_D = D(-\Delta_v)^{-k/2}$ is
of weak type $(1, 1)$ if and only if  $q \leq 2$.
When $q \geq 3$, there exists a constant $C = C(v, D)$
 such that
for all $f \in L(1 + \ln^+L)^{\frac{q}{2} -1}(\mu)$
and all $\lambda > 0$, we have
{\em\begin{align} \label{rf1}
\mu\left\{ x; \left|R_D f(x) \right| > \lambda \right\}
\leq C \int \frac{|f|}{\lambda}
\left( 1 + \ln^+{\frac{|f|}{\lambda}} \right)^{\frac{q}{2} - 1} \, d\mu.
\end{align}}
This inequality is sharp in the sense that $q$ cannot be replaced by any
smaller number.
\end{theo}

Let $(e^{t \Delta_{v} })_{t > 0}$ denote the heat semigroup on
$\R^{(n, v)}$, which is a symmetric diffusion semigroup in the sense
of \cite{St70} (the conservation property can be justified by
Theorem~11.8 in \cite{G09},  and the other properties are obvious).
Further, let $(e^{-t \sqrt{-\Delta_v}})_{t > 0}$ denote the Poisson
semigroup. For $f \in C_0^{\infty}$, we define the
vertical Littlewood-Paley-Stein functions associated with the operator $D$
from \eqref{D} as
\begin{align*}
\mathcal{H}_{D}(f)(x)  &= \left( \int_0^{+\infty}
\left| t^{\frac{k}{2}} D
e^{t \Delta_v} f(x) \right|^2 \frac{dt}{t} \right)^{\frac{1}{2}}
\end{align*}
and
\begin{align*}
\mathcal{G}_{D}(f)(x)  &= \left( \int_0^{+\infty} \left| t^k  D e^{-t \sqrt{-\Delta_v}} f(x) \right|^2 \frac{dt}{t} \right)^{\frac{1}{2}}.
\end{align*}
 The $L^p$ ($1 < p < +\infty$)
boundedness of
$\mathcal{G}_D$ and $\mathcal{H}_D$ is easy to verify.

\begin{theo}\mbox{\textrm{\textbf{(a)}}}\label{H} \hskip3pt
  The operator $\mathcal{H}_D$
is of weak type $(1, 1)$ if and only if $q \leq 1$. When  $q > 1$,
 there exists a constant $C = C(v, D)$ such that for all $f \in L(1 + \ln^+L)^{\frac{q}{2}- \frac{3}{4}}(\mu)$ and all $\lambda > 0$
\begin{align}
\mu\left\{ x; \mathcal{H}_D f(x) > \lambda \right\}
\leq C \int \frac{|f|}{\lambda} \left( 1 +
\ln^+{\frac{|f|}{\lambda}} \right)^{{ \frac{q}{2}-  \frac{3}{4}}} \, d\mu. \label{rf5'}
\end{align}
\mbox{\textrm{\textbf{(b)}}} \hskip3pt The operator $\mathcal{G}_D$
is of weak type $(1, 1)$ if and only if $q \leq 2$. When  $q > 2$,
 there exists a constant $C = C(v, D)$ such that for all
$f \in L(1 + \ln^+L)^{\frac{q}{2} - 1}(\mu)$ and all
$\lambda > 0$
\begin{align}
\mu\left\{ x; \mathcal{G}_D f(x) > \lambda \right\}
&\leq C \int \frac{|f|}{\lambda} \left( 1 +
\ln^+{\frac{|f|}{\lambda}} \right)^{\frac{q}{2} - 1} \, d\mu. \label{rf4}
\end{align}
In \eqref{rf5'} and  \eqref{rf4},
 $q$ cannot be replaced by any smaller number.
\end{theo}

Consider now the horizontal
 Littlewood-Paley-Stein functions and related maximal operators,
  defined for $f \in C_0^{\infty}$ by
\begin{align*}
h_k(f)(x) &= \left( \int_0^{+\infty} \left| t^k
\frac{\partial^k}{\partial t^k} e^{t \Delta_v} f(x) \right|^2
\frac{dt}{t} \right)^{\frac{1}{2}}, \qquad  \quad H_k(f)(x) = \sup_{t > 0}
\left| t^k \frac{\partial^k}{\partial t^k} e^{t \Delta_v} f(x)
\right|,\\
g_k(f)(x) &= \left( \int_0^{+\infty} \left| t^k \frac{\partial^k}{\partial t^k} e^{-t \sqrt{-\Delta_v}} f(x) \right|^2 \frac{dt}{t} \right)^{\frac{1}{2}}, \quad G_k(f)(x) = \sup_{t > 0} \left| t^k \frac{\partial^k}{\partial t^k} e^{-t \sqrt{-\Delta_v}} f(x) \right|.
\end{align*}
Here $k\ge 1$ for $h_k$ and  $g_k$, but  $k\ge 0$ for $H_k$ and
$G_k$.

These operators are bounded on $L^p$ for $1 < p <
+\infty$ in the setting of a general symmetric semigroup with the
contraction property, see \cite{St70}.
The weak type $(1, 1)$ property of $G_k,\; k\ge 0$, in a general setting is
obtained in \cite{LS15}.
For the other three operators,
 we have the following endpoint
estimates. We remark that
 the weak type  $(1,1)$ of  $H_0$ was obtained in
\cite[Theorem 2]{LSW14}.

\noindent \begin{theo} \label{TH2} The operators
  $h_1$,   $H_1$ and
$g_k$ with $k \geq 1$
 are of weak type
$(1, 1)$. For $k \geq 2$, $h_k$ and $H_k$ are
not of weak type $(1,1)$;
 however, there exists a constant $C = C(v, k)$ such that
  for all
$f \in L(1 + \ln^+L)^{ \frac{k}{2}- \frac{3}{4}}(\mu)$
\begin{align}
\mu\left\{ x; h_k f(x)  > \lambda \right\} &\leq C \int
\frac{|f|}{\lambda} \left( 1 + \ln^+{\frac{|f|}{\lambda}}
\right)^{\frac{k}{2}- \frac{3}{4} } \, d\mu
\qquad \forall \lambda > 0, \label{rf2}
\end{align}
and for all $f \in L(1 + \ln^+L)^{\frac{k}{2}-\frac{1}{2}}(\mu)$
\begin{align}
\mu\left\{ x; H_k f(x) > \lambda \right\} &\leq C \int
\frac{|f|}{\lambda} \left( 1 + \ln^+{\frac{|f|}{\lambda}}
\right)^{\frac{k}{2}-\frac{1}{2}} \, d\mu
\qquad  \forall \lambda > 0. \label{rf3}
\end{align}
In these two estimates, the exponents cannot be replaced by any smaller numbers.
\end{theo}

 Our estimates, in particular those involving Orlicz spaces,
are optimal and go beyond earlier known results.

 The structure of this paper is as follows.
Section \ref{sec2} contains estimates for the kernels of the  Riesz transforms,
which are used in Section \ref{sec3} to prove Theorem \ref{one'}. A fundamental
tool here is Proposition \ref{prop2} that is also needed in the later sections.
The two parts of Theorem~\ref{H} are proved in Sections \ref{sec4} and  \ref{sec5}.
Section \ref{sec6} gives estimates for the time derivatives of the heat kernel,
which are applied in the proof of Theorem \ref{TH2} in Section \ref{sec7}.

\medskip

\subsection {Notation and simple facts}

We assume that  $v = e_1 = (1,0,\dots,0)$,
which is no restriction, see \cite{LSW14}.
Then $d\mu = e^{2x_1}dx$,
and it will be convenient to write points in
$\R^n$ as $x = (x_1, x') \in \R \times \R^{n - 1}$.
 Further, we denote $\N =\{0,1,2,\dots \}$
and $\N^+ =\{1,2,\dots \}$.

In the following, $C$  denotes various constants which depend only
 on  $n$ and  $D$ in Sections \ref{sec2} --\ref{sec5},
and only on   $n$ and  $k$ in Sections \ref{sec6} --\ref{sec7}.
 By $A\lesssim B$,
we  mean $A\le C B$ with such a
 $C$ (we say that $A$ is controlled by $B$), and $A\sim B$ stands for  $A\le C B$  and  $B\le C A$.

We have for any $x \in \R^{n}$ and $r \le 1$,
\begin{equation} \label{ESV1}
\mu(B(x, r)) \sim r^n  e^{2 x_1}.
\end{equation}
It follows that our space 
has the local doubling property,

The heat kernel $p_t(x,y)$ in our setting is defined through
\begin{align*}
e^{t \Delta_{e_1}}f(x) = \int  p_t(x,y) f(y)\,d\mu(y), \qquad t>0,   \quad x \in \R^{n},
\end{align*}
for suitable functions $f$. It is explicitly given by (cf.\  \cite{G09} p. 258)
\begin{align}   
p_t(x, y) &= (4 \pi t)^{-\frac{n}{2}} e^{- x_1 - y_1} e^{-t}
e^{-\frac{| x - y |^2}{4 t}}  \label{KLF1}
\\&= (4 \pi t)^{-\frac{n}{2}}
e^{-x_1 - y_1 - |x - y|} \,e^{-t \left(\frac{|x - y|}{2 t} - 1
\right)^2}, \label{KLF2}
\end{align}
for all $t > 0$ and $x = (x_1, x'), \;y = (y_1, y') \in \R \times \R^{n
- 1}$.

Our results include the simpler case $n=1$, but the proofs are written
for  $n\ge 2$. Except for Proposition  \ref{prop2}, we leave it to the
reader to see how the arguments simplify for   $n=1$.

\medskip

\renewcommand{\theequation}{\thesection.\arabic{equation}}
\section{Estimates of the Riesz kernels} \label{sec2}
\setcounter{equation}{0}

\medskip

With $D$, \hskip3pt $k$ and  $q$ as in Theorem \ref{one'},
we write $R_D(x,y)$ for the kernel of
the  Riesz transform
$R_D = D (-\Delta_{e_1})^{-k/2}$.
 The multiindex $\alpha$ will be split  as
 $\alpha = (\alpha_1, \alpha')$.

We state and prove local and  global estimates for  $R_D(x,y)$.

 \begin{propo}\label{local1}
The  kernel $R_D(x,y)$ satisfies  the local standard estimates
\begin{equation}
  \label{stest1}
 |R_D(x,y)| \lesssim \frac{1}{\mu(B(x, | x - y |))},
\qquad 0 < | x - y | \le 1.
\end{equation}
and
\begin{equation}  \label{stest2}
 |\nabla_{y}\, R_D(x,y)| \lesssim    \frac{1}{|x-y|\,\mu(B(x, | x - y |))},
\qquad 0 < | x - y | \le 1.
\end{equation}
 \end{propo}

\begin{propo} \label{glob2}
 The  kernel $R_D(x,y)$ satisfies
 \begin{equation} \label{glob1}
    |R_D(x,y)| \lesssim  
e^{-x_1-y_1 -|x-y|}\, |x-y|^{\frac{q-n-1}2}
\left[1 + \left(\frac{|x'-y'|^2}{|x-y|}  \right)^{\frac{k}2}\right],
\qquad | x - y | > 1.
 \end{equation}
 This estimate is sharp in the sense that
there exists a  ball $B\subset \R^{n-1}$, depending only on $n$ and $D$,  such that
\begin{equation}\label{sharp1}
(-1)^k R_D(x,y) \sim e^{-x_1-y_1 -|x-y|}    \,(x_1-y_1)^{\frac{q-n-1}2}
\qquad \mathrm{for}\qquad \frac{x'-y'}{(x_1-y_1)^{1/2}} \in B,
\end{equation}
if $x_1-y_1  > 0$ is large enough.
\end{propo}

We will need the following integral estimate.

\begin{lem}\label{intest}
  With $\nu \in \R$ and  $a>0$,  let
\[
B_\nu(a) =  \int_0^{+\infty} t^{\nu}\, e^{- t}\, e^{-\frac{a^2}{4 t}} \, \frac{dt}t.
\]
(a) Then for $a\le 2$
\begin{equation*}
 B_\nu(a)  \sim
 \begin{cases}
   a^{2\nu } \quad   &\text{if $\nu < 0$;}\\
     \log (1 + a^{-1}) \quad  &\text{if $\nu = 0$;}\\
    1 \quad  &\text{if $\nu > 0$.}
 \end{cases}
\end{equation*}
(b) Moreover,              
\begin{equation*}
 B_\nu(a) =  \sqrt{2\pi} \,  2^{- \nu} \,  a^{\nu-1/2}\, e^{-a}\,(1+O(a^{-1})), \qquad a \to \infty,
\end{equation*}
and this estimate remains true if  $ B_\nu(a)$ is replaced by
\[
\int_{(a/2-\ell)\vee 0}^{a/2+\ell} t^{\nu}\, e^{- t}\, e^{-\frac{a^2}{4 t}} \, \frac{dt}t.
\]
 for any $\ell \ge a^{3/4}$.

   In this lemma,  the implicit constants  depend on $\nu$ but not on $a$.
\end{lem}

Lemma \ref{intest}  can be proved either with  elementary estimates  or by
evaluating the integral in
terms of Bessel functions according to the 12th equality of \S 3.471
in \cite{GR07}. The latter option also requires an easy estimate of those
parts of the integral left out in the last  statement.

 We will also use a trivial extension of this lemma. Let
$Q(t) = \sum_\nu q_\nu \,t^\nu$ be any finite linear combination of real powers of
the variable $t>0$. Then  Lemma \ref{intest} (b) implies that
\begin{equation}\label{ext}
   \int_0^{+\infty} Q(t) \, e^{- t}\, e^{-\frac{a^2}{4 t}} \, \frac{dt}t
 =  \sqrt{2\pi} \:  Q(a/2)\, a^{-1/2}\, e^{-a}\,(1+O(a^{-1})), \qquad a \to \infty.
\end{equation}

To prepare for the proof of the two propositions, we deduce an expression for
 $R_D(x,y)$. From \eqref{KLF1} it follows that
the kernel of $(-\Delta_{e_1})^{-k/2}$ is
\begin{align*}
  \label{fractker}
 (-\Delta_{e_1})^{-\frac k2}(x,y)   =&\, \frac1{\Gamma({ k/2})}\, \int_0^{+\infty} t^{\frac{k}2}\,
p_t(x,y)\, \frac{dt}t   \\=& \, (4\pi)^{-\frac n2}\,
\frac1{\Gamma({ k/2})}\,    e^{-x_1-y_1}
\int_0^{+\infty} t^{\frac{k-n}2} e^{- t} e^{-\frac{|x-y|^2}{4 t}} \, \frac{dt}t.
\end{align*}
Further,
\begin{equation}\label{RR}
R_D(x,y) = D_x (-\Delta_{e_1})^{-\frac k2}(x,y) = \frac1{\Gamma({ k/2})}\, \int_0^{+\infty} t^{\frac{k}2}\,
D_x p_t(x,y)\, \frac{dt}t.
\end{equation}

 Recall that the Hermite polynomial of degree $j \in \N$  is defined by
\begin{align*}
H_j(s) = (-1)^j \,e^{s^2}\, \frac{d^j}{ds^j} \, e^{-s^2},
\end{align*}
and $H_\alpha = \otimes_{i}\, H_{\alpha_i}$ for any multiindex $\alpha$.
It is well known that the $H_\alpha$ are
orthogonal with respect to the Gaussian measure $e^{-|x|^2}dx$ in $\R^n$ and that the
 leading term of   $H_j(s)$ is $2^js^j$.

 The definition of  $H_\alpha$ implies that
\begin{align}\label{deriv}
 \partial^{\alpha}_x\, e^{-\frac{|x-y|^2}{4 t}} =
\left(-\frac{1}{2 \sqrt{t}}\right)^{|\alpha|}\,
H_\alpha\left(\frac{x - y}{2 \sqrt{t}}\right) e^{-\frac{|x-y|^2}{4 t}}.
\end{align}

The kernel  $R_D(x,y)$ is a linear combination of terms
$\partial^{\alpha}_x\,(-\Delta_{e_1})^{-k/2}(x,y) $, with
$\alpha = (\alpha_1,\alpha')$ a multiindex of length  $|\alpha| = k$
and with $\alpha_1 \le q$.
When we differentiate the expression in \eqref{KLF1},
some differentiations with respect to $x_1$ may fall on the factor
$e^{-x_1}$. Those differentiations falling on
the factor $e^{-\frac{|x-y|^2}{4 t}}$ will be given by a multiindex
$\tilde \alpha = ( \tilde \alpha_1,\alpha') $, with
$0 \le \tilde \alpha_1 \le \alpha_1.$
We then see from \eqref{deriv}, with $\alpha$ replaced by $\tilde \alpha$,
that $\partial^{\alpha}_x\,(-\Delta_{e_1})^{-k/2}(x,y) $  will be a sum of terms
\begin{align}\label{term7}
&(-1)^{\alpha_1 - \tilde \alpha_1}\,  e^{-x_1-y_1} \int_0^{+\infty} t^{\frac{k-n}2} e^{- t}
\left(-\frac{1}{2 \sqrt{t}}\right)^{|\tilde \alpha|}\,
H_{\tilde \alpha}\left(\frac{x - y}{2 \sqrt{t}}\right) e^{-\frac{|x-y|^2}{4 t}}
\,\frac{dt}t \notag \\ =\,&
(-1)^{|\alpha|}\, 2^{-|\tilde \alpha|}\, e^{-x_1-y_1} \int_0^{+\infty}
t^{\frac{\alpha_1 -\tilde \alpha_1-n}2} H_{\tilde \alpha_1}\left(\frac{x_1 - y_1}{2 \sqrt{t}}
\right) H_{ \alpha'}\left(\frac{x' - y'}{2 \sqrt{t}}\right) e^{- t}\,
e^{-\frac{|x-y|^2}{4 t}} \,\frac{dt}t,
\end{align}
where $\tilde \alpha_1$ runs from 0 to  $ \alpha_1$ and the terms have
positive coefficients coming from Leibniz' formula
and from the factors $(4\pi)^{-n/2}$ and $1/\Gamma(k/2)$.

Expanding the Hermite polynomials here, we obtain a sum of terms proportional
to
\begin{align}\label{term}
 e^{-x_1-y_1}  (x_1 - y_1)^{\gamma_1}  (x' - y')^{\gamma'}  \int_0^{+\infty}
t^{\frac{\alpha_1 -\tilde \alpha_1-\gamma_1-|\gamma'|-n}2}  e^{- t}\, e^{-\frac{|x-y|^2}{4 t}}
\,\frac{dt}t,
\end{align}
and the sum is now  taken also over a  multiindex  $(\gamma_1,\gamma') $
 with $0 \le \gamma_1   \le \tilde \alpha_1 $
and
$0 \le \gamma'   \le  \alpha' $ (componentwise ordering).

\vskip4pt

\noindent {\bf Proof of Proposition \ref{local1}.}
  Here $|x-y| \le 1$.
If the exponent of $t$
in the  integral in \eqref{term} is negative, we see from Lemma \ref{intest} (a) that the modulus
of  the   expression  \eqref{term} is controlled by
\begin{align*}
 e^{-2x_1}  |x-y|^{\gamma_1+|\gamma'|}
|x-y|^{\alpha_1 -\tilde \alpha_1-\gamma_1-|\gamma'|-n}
\le e^{-2x_1}  |x-y|^{-n}.
\end{align*}
For other values of the exponent, the bound $ e^{-2x_1}|x-y|^{-n}$
also follows. This is the first standard estimate \eqref{stest1}. To obtain
also \eqref{stest2}, it is enough to trace the argument given, with a
 differentiation also in $y$.
\hfill $\Box$

\vskip4pt

\noindent {\bf Proof of Proposition \ref{glob2}.}
For  $|x-y| > 1$,   Lemma \ref{intest} (b)
implies that the expression \eqref{term} equals
 constant times
\begin{align}\label{term1}
 e^{-x_1-y_1-|x-y|}  (x_1-y_1)^{\gamma_1}  (x' - y')^{\gamma'}
|x-y|^{\frac{\alpha_1 -\tilde \alpha_1-\gamma_1-|\gamma'|-n-1}2}\,(1+O(|x-y|^{-1})),
\end{align}
whose modulus is controlled by
\begin{align}\label{term6}
 e^{-x_1-y_1-|x-y|} \,
|x-y|^{\frac{\alpha_1 -(\tilde \alpha_1-\gamma_1) -n-1}2}  \left(\frac{|x' - y'|^2}{|x-y|}  \right)^{\frac{|\gamma'|}2}.
\end{align}
   Now \eqref{glob1} follows,
because $\alpha_1 \le q$ and $ \gamma_1   \le \tilde \alpha_1 $
and also $0 \le |\gamma'|   \le  |\alpha'| = k-\alpha_1 $.

To verify  \eqref{sharp1}, assume that
$x_1-y_1 > 0$ is large and that $|x' - y'|^2/(x_1-y_1)$ stays bounded. In this
argument, we will neglect all terms which are much smaller
than the right-hand side of \eqref{sharp1}. From  \eqref{term1} and
\eqref{term6}, we then see that   in  \eqref{term7} we  need only
take  $\alpha_1 = q$, which implies  $|\alpha'| = k-q$, and $ \gamma_1   = \tilde \alpha_1 $.
The latter equality means that in
the Hermite polynomial $H_{\tilde \alpha_1}$ in \eqref{term7}, we
  consider only the leading term, which is
 $$2^{\tilde \alpha_1}\,\left(\frac{x_1 - y_1}{2 \sqrt{t}}\right)^{\tilde \alpha_1}
= t^{-\frac{\tilde \alpha_1}2} \,(x_1 - y_1)^{\tilde \alpha_1}.$$


Instead of \eqref{term7} we will now have
\begin{align*}     
 (-1)^{k}\, 2^{-|\tilde \alpha|} e^{-x_1-y_1}   (x_1 - y_1)^{\tilde \alpha_1}
 \int_0^{+\infty} t^{\frac{q -2\tilde \alpha_1-n}2}
H_{\alpha'}\left(\frac{ x' - y'}{2\sqrt t} \right)
  e^{- t}\, e^{-\frac{|x-y|^2}{4 t}} \,\frac{dt}t.
\end{align*}
  Applying  \eqref{ext} to the integral here, we see that
\eqref{term7} amounts to a  positive constant times
\begin{align}  \label{term5}
 (-1)^{k}\, 2^{-|\tilde \alpha|}\ e^{-x_1-y_1-|x-y|} \, \notag
  (x_1 - y_1)^{\tilde \alpha_1}\,
\left(\frac{|x-y|}2\right)^{ \frac{q - 2 \tilde \alpha_1 - n}{2}}
 & |x-y|^{-\frac12}\,H_{\alpha'}\left(\frac{ x' - y'}{\sqrt{2}\, |x-y|^{1/2} } \right)
 \\ &\times  (1+O(|x-y|^{-1})).
\end{align}
Here we can replace the powers of $|x-y|$ by the same powers of
$x_1-y_1$, since
 $|x-y| = x_1-y_1 +O(1)$, so \eqref{term5} equals
\begin{align*}   
 (-1)^{k}\,  2^{- \frac{q}{2}+ \tilde \alpha_1 -|\tilde \alpha|  + \frac{n}{2}} \, e^{-x_1-y_1-|x-y|}\, (x_1-y_1)^{\frac{q-n-1}2}
H_{ \alpha'}\left(\frac{x' - y'}{{ \sqrt{2}}\, (x_1-y_1)^{1/2}}\right)
  (1+O(|x-y|^{-1})).
\end{align*}

Summing over $0\le\tilde \alpha_1\le q $, we conclude that
\begin{multline*}   
 \partial^{\alpha}_x\,(-\Delta_{e_1})^{-k/2}(x,y)\\ = (-1)^{k}\,  b_{\alpha}\, e^{-x_1-y_1-|x-y|}\, (x_1-y_1)^{\frac{q-n-1}2}
H_{ \alpha'}\left(\frac{x' - y'}{{ \sqrt{2}}\, (x_1-y_1)^{1/2}}\right)
 (1+O(|x-y|^{-1}))
\end{multline*}
for some $b_{\alpha}>0$.
Here $\alpha = (q, \alpha')$ and  $|\alpha| = k$, and if  we sum over such
$\alpha$ with the coefficients from  \eqref{D}, the result will be
\begin{align*}   
R_D(x,y)&\\ =\,  & (-1)^{k} e^{-x_1-y_1-|x-y|} (x_1-y_1)^{\frac{q-n-1}2}
\sum a_{\alpha}  b_{\alpha} H_{ \alpha'}\left(\frac{x' - y'}{{ \sqrt{2}}\, (x_1-y_1)^{1/2}}\right)
+ \; \mathrm{negligible}\; \mathrm{terms.}
\end{align*}
Since the $a_{\alpha}$ do not all vanish, the orthogonality property
of the Hermite
polynomials implies that the polynomial given by the sum here is not
identically 0. To finish the proof of   \eqref{sharp1}, we need only take
a closed ball  $B$ where this polynomial does not vanish.
\hfill $\Box$

\vskip19pt

\renewcommand{\theequation}{\thesection.\arabic{equation}}
\section{Proof of Theorem \ref{one'} } \label{sec3}
\setcounter{equation}{0}

\medskip

We split $R_D$ into a  part at infinity
\begin{equation*}
   R_D^{\infty}f(x) = \int_{|x-y| > 1}R_D(x,y) f(y)\,d\mu(y)
\end{equation*}
 and a local part $R_D^{\mathrm{loc}} = R_D -R_D^{\infty}$.

The local part is easy to treat.
Because of the  local doubling property,
  we can use the  method of localization. In view of Proposition \ref{local1},
standard Calder\'on-Zygmund singular integral  theory gives the weak
type $(1, 1)$ of $R_D^{\mathrm{loc}}$.

To estimate $R_D^{\infty}$, we start with that part defined by the
restriction  $x_1-y_1 \le 1$. Then we have integrability, since
  \eqref{glob1} implies
\begin{multline*}    
\int_{\substack{|x-y| > 1\\x_1-y_1 \le 1}}    
|R_D(x,y)|\,d\mu(x)  \lesssim
\int_{|x-y| > 1} e^{1-|x-y|}\, |x-y|^{\frac{q -n-1}2}
\left[1 + \left(\frac{|x'-y'|^2}{|x-y|}  \right)^{\frac{k}2}\right]
\,dx \lesssim 1,
\end{multline*}
uniformly in $y$. It follows that this part of the operator is of strong type $(1,1)$.

To treat the opposite case  $x_1-y_1 > 1$, we write the first exponent in  \eqref{glob1} as
\[
-2x_1-(|x-y|-(x_1-y_1))
\]
 and observe that in this case
\begin{align*}
 |x-y|-(x_1-y_1) = \frac{|x'-y'|^2}{|x-y|+x_1-y_1}
\end{align*}
and     
\begin{align} \label{compar}
 \frac{|x'-y'|^2}{2|x-y|} \le  \frac{|x'-y'|^2}{|x-y|+x_1-y_1} \le  \frac{|x'-y'|^2}{|x-y|}.
\end{align}

Thus Proposition \ref{glob2} implies that for  $x_1- y_1 > 1$
\begin{align*}
|R_D(x,y)| \lesssim e^{-2x_1}\, |x-y|^{\frac{q-n-1}2}
\exp\left( - \frac{|x'-y'|^2}{2|x-y|}  \right)
\left[1 + \left(\frac{|x'-y'|^2}{|x-y|}  \right)^{\frac{k}2}\right].
\end{align*}
It follows that this part of $|R_D(x,y)|$  is controlled by the kernel
\begin{align}\label{defM}
\mathcal V_{\kappa}(x,y) =  e^{-2x_1} \,|x-y|^{\frac{\kappa-n -1 }2}
\exp\left( -\frac14 \frac{|x'-y'|^2}{|x-y|}  \right)\,\chi_{\{x_1- y_1 > 1\}},
\end{align}
with $\kappa = q.$

\vskip15pt

Let for  $\kappa \in \R$ and suitable functions $f$
\begin{align}
\widetilde{\mathcal{V}_\kappa}f(x) = \int_{\R^{n}} \mathcal{V}_\kappa(x, y)\, f(y) \,d\mu(y). \label{DT2}
\end{align}

To prove Theorem \ref{one'} in the case  $q \leq 2$,
it is enough to show the following result.

\begin{propo} \label{prop2}
The operator $\widetilde{\mathcal{V}_2}$ defined by \eqref{DT2} is of weak type
$(1, 1)$.
\end{propo}

 The case $n = 1$ is trivial, since then
\begin{align*}
|\widetilde{\mathcal{V}_2}f(x)| \leq e^{-2 x} \| f \|_{L^1(d\mu)},
\end{align*}
which  implies the weak type $(1, 1)$ of
$\widetilde{\mathcal{V}_2}$.

In what follows, we suppose that $n \geq 2$, so that
\[
\widetilde {\mathcal{V}_2}f(x) =e^{-2x_1} \int_{y_1 < x_1 - 1}   |x-y|^{\frac{1-n}2}\,
\exp\left(-\frac  {1}{4}\, \frac {|x'-y'|^2 }{|x-y|}
\right)\, f(y)e^{2y_1}\,dy,
\]
where we may assume $0\le f \in L^{1}(d\mu) $.  It will be convenient
to use $g(y) = f(y)e^{2y_1} \in L^1(dy)$ instead of $f$.
When $|x'-y'| \ls \sqrt{x_1 - y_1}$, the exponential in the
integrand  is essentially 1. This indicates the most important
part of the operator, dealt with in  the following proposition.

\begin{propo}\label{pr8}
The operator
\begin{equation}\label{prop8}
Tg(x) = e^{-2x_1} \int_{\substack{y_1 < x_1 - 1\\|x'-y'| <  
\sqrt{x_1 - y_1}}}\,
  (x_1-y_1)^{\frac{1-n}2} g(y)\,dy
\end{equation}
maps  $L^1(dy)$ boundedly into  $L^{1,\infty}(d\mu)$.
\end{propo}

We first verify that this proposition implies Proposition \ref{prop2}. Let for $j = 1, 2, \dots$
\begin{multline*}
T_jg(x) = e^{-2x_1}
\int_{\substack{y_1 < x_1 - 1\\  2^{j-1} \sqrt{x_1 - y_1} \le |x'-y'| < 2^{j} \sqrt{x_1 - y_1}}}\,
  (x_1-y_1)^{\frac{1-n}2}\,\exp\left(-\frac  {1}{4}\,
\frac {|x'-y'|^2 }{|x-y|}  \right)\,  g(y)\,dy;
\end{multline*}
then $\widetilde{\mathcal{V}_2}f \le Tg + \sum_1^\infty T_jg $. In this integral, the
exponential is less than $\exp(-c2^j)$ for some $c>0$, so that
\[
T_jg(x) \le \exp(-c2^j)\, e^{-2x_1}  \int_{\substack{y_1 < x_1 - 1 \\|x'-y'|
< 2^j \sqrt{x_1 - y_1}}}\,
  (x_1-y_1)^{\frac{1-n}2}\, g(y)\,dy.
\]
Scaling the variables $x'$ and $y'$ by a factor $2^j$, we can
control this expression in terms of the operator $T$, still with a
rapidly decreasing coefficient. One can then apply  Proposition~\ref{pr8}
and use \cite[Lemma 2.3]{SW69} to sum in weak $L^1$. The result will
be the weak type $(1,1)$  of $\widetilde{\mathcal{V}_2}$, which is
Proposition \ref{prop2}.

   \vskip4pt

\noindent \textbf{Proof of Proposition \ref{pr8}}.
 \hskip3pt We cover $\rr^{n-1}$ and  $\rn$ with  lattices of unit cubes
\[
Q_{ m} = \left\{x' = (x_2,\dots,x_n) \in \rr^{n-1}: m_\nu - 1 < x_\nu \le m_\nu ,\;
\; \nu=2,\dots,n \right\},    
\]
where $m = (m_2, \dots, m_n) \in \zz^{n-1}$,
and
\[
Q_{i, m} =  (i-1, i] \times Q_{ m}, \qquad (i,m)\in \zz \times \zz^{n-1}.
\]
In this proof, the word ``cube'' will always refer to cubes of
these lattices.

 Given $0\le g  \in L^1(dy)$ and $\lambda > 0$, we must
estimate the $\mu$ measure of the level set $L_\lambda = \{x: Tg(x) >\lambda  \}$. Since
$ Tg(x) \le  e^{-2x_1} \int g(y)\,dy$, the coordinate $x_1$ is
bounded above on
 $L_\lambda$. Without restriction, we can then assume that  $L_\lambda$ is contained in the left half-plane
$\{x: x_1 \le 0 \}$.

The aim of the following lemma is to   discretize the  variable $x_1$.

  \begin{lem}\label{hej}
Let  $i \le 0$ and  assume that the cube $Q_{i, m}$ intersects
the level set $L_\lambda$. Then
$\tilde  Tg(i,x') >e^{-2}\,2^{-(n-1)/2}\,\lambda  $
for all $x' \in  Q_{m}$,
 where
\begin{equation}  \label{X}
\tilde Tg(i,x') =  e^{-2i} \int_{\substack{y_1 < i-1 \\ |y'-x'| < \sqrt{i-y_1}
+\sqrt{n-1}}}
   (i-y_1)^{\frac{1-n}2}\,g(y)\,dy.
\end{equation}
  \end{lem}

 \noindent \textbf{Proof.}
Take a point
$z \in Q_{i, m}\cap L_\lambda$, so that
\begin{equation}\label{zz}
 \lambda < Tg(z) = e^{-2z_1} \int_{\substack{y_1 < z_1 - 1\\|y'-z'| <  
\sqrt{z_1 - y_1}}}\,
  (z_1-y_1)^{\frac{1-n}2} g(y)\,dy.
\end{equation}
Here  $i < 1+z_1$ and $z_1-y_1 > 1$, so that $i-y_1 < 1+z_1-y_1 < 2(z_1-y_1)$  and thus $(z_1-y_1)^{(1-n)/2} < 2^{(n-1)/2}(i-y_1)^{(1-n)/2}$.
  Further, the region of integration in  \eqref{zz} is contained in that in
\eqref{X}, because $z_1 \le i$ and $|z'-x'|\le \sqrt{n-1}$ for any
$x' \in  Q_{m}$. Since also $e^{-2z_1} \le e^2 \, e^{-2i}$,
it follows that $Tg(z) <  e^2\, 2^{(n-1)/2}\, \tilde Tg(i,x')$,
which proves the lemma. ~ \hspace*{20pt} ~ \hfill $\Box$

\vskip4pt

The region of integration in \eqref{X} is contained in the set
\[
D_{i, m} =  \bigcup_{x'\in Q_m}
\left\{(y_1, y'): y_1 < i-1,\; |y'-x'| < \sqrt{i - y_1}+ \sqrt{n-1}\right\}.
\]

For each cube $Q_{i, m}$, we define a family of ``forbidden'' cubes
\[
F_{i, m} =
\left\{Q_{i', m'}: i' < i,\; |m-m'| < 4\sqrt{i - i'}+ 4\sqrt{n-1}\right\}
\]
whose union contains $D_{i, m}$, that is,
\[
D_{i, m} \subset \hat D_{i, m} :=  \bigcup_{Q\in F_{i, m}}  Q,
\]
as easily verified.
The exponential behavior of  $\mu$ implies that
 \begin{equation}
   \label{forb}
\mu(\hat D_{i, m}) \lesssim \mu( Q_{i, m}).
 \end{equation}

For each $i\le 0$, we will select some of the cubes  $Q_{i, m}$, namely
those for which $m$ is in a set $A_i\subset \zz^{n-1}$.
 These $A_i$ will be
 defined by recursion. We will have the inclusion
 \begin{equation}
   \label{I}
   L_\lambda   \subset
\bigcup_{i\le 0} \:  \bigcup_{m \in A_i} \left(Q_{i, m} \cup \hat D_{i, m}\right).
 \end{equation}
Because of \eqref{forb}, this implies that
\begin{equation}
  \label{enough}
\mu(L_\lambda) \lesssim \sum_{i\le 0}  \sum_{m \in A_i} \mu( Q_{i, m})
\sim \sum_{i\le 0}  \sum_{m \in A_i} e^{2i}.
\end{equation}

When we  make the selection, i.e.,  define the  $A_i$, we will consider
the $\hat D_{i, m}$ as forbidden regions, where no selection is allowed.
Indeed, we will select  for each $i = 0, -1, -2, \dots$ those cubes $Q_{i, m}$
which intersect the level set and are not in $F_{j, \ell}$
for any  already selected cube $ Q_{j, \ell}$. In this way, cubes selected
at different steps will be far from each other.

More precisely, we first let
\[
A_0 =  \left\{ m \in \zz^{n-1}:  Q_{0, m} \cap L_\lambda
\ne \emptyset\right\}.
\]
Assuming  for some $i<0$ the sets $A_{i+1}, \dots A_0$ already defined,
we then let
\[
A_i =  \left\{ m \in \zz^{n-1}:  Q_{i, m} \cap L_\lambda
\ne \emptyset   \; \; \; \mathrm{and}
\; \; \; Q_{i, m} \notin F_{j, \ell}  \; \; \; \mathrm{for}  \; \; \mathrm{all}
\; \; \; j> i, \; \; \ell \in A_{j}  \right\}.
\]
The inclusion \eqref{I} is immediate from this construction.

 For each $i\le 0$, we define
  \[
D_i = \bigcup_{m \in A_i}  D_{i, m},
\]
which is a union of cubes contained in $\{y: y_1 \le i-1 \}$.
In order to replace the $D_i$ by pairwise disjoint sets, we let $E_0
= D_0$ and
\[
E_i =  D_i\setminus \bigcup_{j =i+1}^0 D_{j}, \qquad i= -1,-2,\dots.
\]
Then we let $g_i = g\, \chi_{E_i}$ and observe that
 $g = \sum_{j=i}^0 g_j$ in $D_i$.           
This together with Lemma \ref{hej} implies that, for
 $i \le 0$ and  $x' \in Q_{m} $ with $m \in A_i$,
 \begin{equation}
   \label{cruc}
    \lambda \lesssim \tilde Tg(i,x')  =
       e^{-2i} \sum_{j=i}^0 \int_{y_1 < i-1} (i-y_1)^{\frac{1-n}2}  \int_{ |y'-x'| < \sqrt{i-y_1}+\sqrt{n-1}}\,
 g_j(y)\,dy'dy_1.
\end{equation}

The integral in $y'$ in this expression can be seen as a convolution in
$\rr^{n-1}$. That leads to an estimate which is insufficient for our purpose,
essentially because several values of $i$ are used for the same $g_j$.
To do better,  we will  replace  $i$ by $j$
in the iterated integral in  \eqref{cruc}. Thus
we first  claim that
\begin{equation}
  \label{repl}
  (i-y_1)^{\frac{1-n}2} g_j(y) \lesssim   (j-y_1)^{\frac{1-n}2} g_j(y)
\end{equation}
if $y,\:x',\:m$, $i$ and $j$ are as in   \eqref{cruc}.
Then we observe that  an immediate
consequence of  \eqref{cruc} and  \eqref{repl} is that if
$x' \in Q_{m} $ for some $m \in A_i$, then
 \begin{equation} \label{cruc'}
      \lambda \lesssim
  e^{-2i} \sum_{j=i}^0 \int_{y_1 < j-1} (j-y_1)^{\frac{1-n}2}  \int_{ |y'-x'| < \sqrt{j-y_1}
+\sqrt{n-1}}\,
 g_j(y)\,dy'dy_1,
\end{equation}
where we also made the regions of integration larger.

Before we use \eqref{cruc'} to finish the proof of the lemma,
 we  verify the claim \eqref{repl}. Then we
can obviously assume that $g_j(y)\ne 0$, which implies $y \in E_j \subset D_j$.
The definition of $D_j$  says that
 $y_1<j-1$ and
$   |y'-z'| < \sqrt{j - y_1}+ \sqrt{n-1}$
for some $z'\in Q_\ell$ and some $\ell\in A_j$. Clearly $|z'-\ell| \le \sqrt{n-1}$,
so the triangle inequality leads to
\begin{align}\label{14}
   |y'-\ell| < \sqrt{j - y_1}+ 2\sqrt{n-1}.
\end{align}
The inner integral in \eqref{cruc} is taken over those $y'$ satisfying
\begin{align}\label{16}
  |y'-x'| < \sqrt{i-y_1}+\sqrt{n-1}.
\end{align}

Since also  $|x'-m| \le\sqrt{n-1} $,  \eqref{14} and \eqref{16} imply
\begin{align}\label{15}
   |m-\ell| < \sqrt{j - y_1}+ \sqrt{i - y_1}+ 4\sqrt{n-1}
\le 2\sqrt{j - y_1}+ 4\sqrt{n-1}.
\end{align}

On the other hand, the cube $Q_{i,m}$ was selected and thus not forbidden
by the selected cube $Q_{j,\ell}$, i.e.,  $Q_{i,m} \notin F_{j,\ell}$. This means
that $|m-\ell| \ge 4\sqrt{j - y_1}+ 4\sqrt{n-1}$. Combining this inequality
with \eqref{15}, we get $4\sqrt{j-i} <  2\sqrt{j-y_1}$, so that
 $j-i <  (j-y_1)/4$.  Then we can write
\[
j-y_1 = j-i+i-y_1 < \frac14 (j-y_1) + i-y_1
\]
and thus
\[
j-y_1  < \frac43 (i-y_1).
\]
This implies  \eqref{repl}, and  \eqref{cruc'} also follows.

We now  integrate the inequality \eqref{cruc'} in $x'$ over the unit cube $Q_m$, to conclude that
 \begin{equation}
   \label{cruc2}
 e^{2i}    \lesssim \frac1 \lambda  \sum_{j=i}^0
\int_{y_1 < j-1} \int_{Q_m} (j-y_1)^{\frac{1-n}2} \int_{ |y'-x'| < \sqrt{j-y_1}
+\sqrt{n-1}}\,
 g_j(y)\,dy'dx'dy_1.
\end{equation}
The inner  part of this expression,
\[
G_j(y_1,x'): =  (j-y_1)^{\frac{1-n}2} \int_{ |y'-x'| < \sqrt{j - y_1}+ \sqrt{n-1}}
g_j(y_1,y')\,dy',
\]
is the value at $x'$ of the convolution in $\rr^{n-1}$ of $g_j(y_1,.)$
and the characteristic
function of the ball
$\{y' \in\rr^{n-1}: |y'| < \sqrt{j - y_1}+ \sqrt{n-1}  \}$,
essentially normalized in  $L^1(\rr^{n-1})$ by the first factor.
Thus
\begin{equation}\label{Gg}
\int_{\rr^{n-1}}  G_j(y_1,x') \,dx' \lesssim \int_{\rr^{n-1}} g_j(y_1,y')\,dy',
\end{equation}
uniformly in $j-y_1 > 1$.

Combining \eqref{cruc2} with \eqref{enough}, we have
 \begin{align*}
 \mu(L_\lambda) \lesssim  \frac1 \lambda \sum_{i\le 0}  \sum_{m \in A_i} \sum_{j=i}^0
\int_{y_1 < j-1} \int_{Q_m} G_j(y_1,x')\, dx'dy_1
\\ = \frac1 \lambda  \sum_{j\le 0}
\int_{y_1 < j-1} \sum_{i = -\infty}^j \sum_{m \in A_i} \int_{Q_m} G_j(y_1,x')\, dx'dy_1.
\end{align*}

In the last expression here, we have a double sum over all $m \in A_i$ and
$i\le j$. The corresponding family of  cubes $Q_m$ is pairwise disjoint, since
  $m \in A_i$ implies $Q_{i',m} \in F_{i,m}$
for $i'<i$ by the definition of  $F_{i,m}$. Thus the double sum of integrals
over $Q_m$ is no larger than the left-hand side of \eqref{Gg}.
As a result,
\begin{align*}
 \mu(L_\lambda) \lesssim  \frac1 \lambda \sum_{j\le 0}
\int_{y_1 < j-1}  \int_{\rr^{n-1}} g_j(y_1,y')\,dy'dy_1  \le
\frac1\lambda \sum_{j\le 0}  \int_{\rr^{n}} g_j(y)\,dy
\le \frac1\lambda \int g(y)\,dy.
\end{align*}
 This ends the proof of Proposition \ref{pr8}.
\hfill $\Box$

\vskip7pt

Consider now the case $q > 2$ in Theorem \ref{one'}.
The main part of the kernel $R_D(x,y)$ is controlled by
the kernel $\mathcal V_q$ defined in \eqref{defM}. Therefore,
the estimate  \eqref{rf1}
is a consequence of the following proposition.

\begin{propo} \label{propn1} For any $\kappa > 2$,
there exists a constant $C = C(\kappa, n) > 0$ such that
for all $f \in L(1 + \ln^+L)^{{\kappa}/{2} - 1}(\mu)$ and $\lambda > 0$, we have
{\em\begin{align} \label{rf1n}
\mu\left\{ x; \left|\widetilde{\mathcal{V}_\kappa} f(x)\right| > \lambda \right\} \leq C \int \frac{|f|}{\lambda}
\left( 1 + \ln^+{\frac{|f|}{\lambda}} \right)^{\frac{\kappa}{2} - 1} \, d\mu,
\end{align}}
where the operator $\widetilde{\mathcal{V}_\kappa}$ is defined by \eqref{DT2}.
\end{propo}

\noindent {\bf Proof.} We have for $h \in L^1(\mu)$
\begin{align*}
\widetilde{\mathcal{V}_\kappa}h(x) = \int_{y_1<x_1-1}\,   e^{-2(x_1-y_1)} \,|x-y|^{\frac{\kappa-n -1 }2}
\exp\left( -\frac14 \frac{|x'-y'|^2}{|x-y|}  \right)\,h(y) \,
dy.    
\end{align*}
Here we assume $h\ge 0$ and apply the elementary inequality
\begin{align*} 
a b^{\frac{\kappa}{2} - 1} \leq C_0\, \left[ a (1 + \ln^+{a})^{\frac{\kappa}{2} - 1} + e^{\frac{b}{8}} \right], \qquad a, b > 0,
\end{align*}
easily proved by separating the cases $b \le 10(1+\ln^+{a})$ and
 $b > 10(1+\ln^+{a})$. Here $C_0$ depends only on $\kappa$. Letting $a = h(y)$ and $b =|x-y|$, we get
\begin{multline*}
\widetilde{\mathcal{V}_\kappa}h(x) \le
C_0\,\widetilde{\mathcal{V}_2}\left(h(1+\ln^+{h})^{\frac{\kappa}{2}-1}\right)(x)\\
+ C_0\,
 \int_{y_1<x_1-1}   e^{-2(x_1-y_1)} \,|x-y|^{\frac{1-n}2}
\exp\left( -\frac14 \frac{|x'-y'|^2}{|x-y|}  \right)\, e^{\frac{|x-y|}8}\,dy.
\end{multline*}
The  integral here equals     
\begin{equation*}        
  \int_{z_1>1}   e^{-2z_1} \,|z|^{\frac{1-n}2}
\exp\left( -\frac14 \frac{|z'|^2}{|z|}  \right)\, e^{\frac{|z|}8}\,dz,
\end{equation*}
and to see that it is finite, we split it in two:
 \begin{equation*}
 \int_{\substack{z_1>1\\ |z'|<2z_1}} + \int_{\substack{z_1>1\\ |z'|>2z_1}} = I_1 + I_2.
\end{equation*}
In $I_1$, we can estimate the second and third exponential factors by $1$ and $e^{3z_1/8}$,
 respectively, and integrate first in $z'$. In  $I_2$, we have $|z| < 3|z'|/2$ and thus
 $|z'|^2/|z| >2|z'|/3$, so that the second exponential in  $I_2$ is less than
$\exp(-|z'|/6)$. It is then enough to estimate the third exponential by
$\exp(z_1/8+|z'|/8)$.

It now follows that
\begin{equation*}
\widetilde{\mathcal{V}_\kappa}h(x) \le
C_0\,
\widetilde{\mathcal{V}_2}\left(h(1+\ln^+{h})^{\frac{\kappa-2}{2}}\right)(x) + C_1
\end{equation*}
for some constant $C_1$.

Given $f \in L(1+\ln^+L)^{{(\kappa-2)}/2}(\mu)$ and $\lambda > 0$, we let
$h = {2C_1}\lambda^{-1} |f|$. Then $|\widetilde{\mathcal{V}_\kappa}f(x)| > \lambda$
implies $\widetilde{\mathcal{V}_\kappa}\left(h\right)(x)> 2C_1$ and
$\widetilde{\mathcal{V}_2}\left(h(1+\ln^+{h})^{(\kappa-2)/{2}}\right)(x) > C_1/C_0$.
Since $\widetilde{\mathcal{V}_2}$ is of weak type $(1,1)$ by Proposition \ref{prop2},
we conclude that
\begin{equation*}
  \mu\{ |\widetilde{\mathcal{V}_\kappa}f(x)| > \lambda\}
\lesssim \left\|h(1+\ln^+{h})^{\frac{\kappa}{2}-1} \right\|_{L^1(\mu)}
\sim \int \frac{|f|}\lambda
\left(1+\ln^+\frac{|f|}\lambda\right)^{\frac{\kappa}2-1}
\,d\mu.
\end{equation*}
Proposition \ref{propn1} is proved.
\hfill $\Box$

\bigskip

To show that  \eqref{rf1} is sharp,  we will use \eqref{sharp1}.

  Let
$f = \chi_{B(0,1)}$. With $\eta > 0$ large,
we  consider $R_Df(x)$ in the region
\begin{align}\label{omeg}
\Omega_\eta = \left\{x: \eta < x_1 <  \eta + 1, \quad
\frac{x'}{\sqrt\eta} \in \frac12\,  B \right\}.
\end{align}
Here  $ \frac12 B$ denotes the concentric scaling of the  ball
$B$ from  \eqref{sharp1},
by a factor $1/2.$
Let $x\in \Omega_\eta$ and $y \in B(0,1)$. Then
\begin{align*}
\frac{x'-y'}{(x_1-y_1)^{1/2}} = \frac{x'}{\sqrt\eta} + O\left(\frac1{\sqrt\eta} \right), \qquad \eta \to \infty.
\end{align*}
This implies that $(x'-y')/(x_1-y_1)^{1/2} \in B$  for large $\eta$,
 and  \eqref{sharp1} applies.
Further, \eqref{compar} shows that $e^{-x_1-y_1-|x-y|} \sim e^{-2x_1} $, and so
\begin{align*}
(-1)^k R_Df(x) \sim e^{-2\eta}\, \eta^{\frac{q-n-1}2} = \lambda
\end{align*}
for  $x\in \Omega_\eta$, where the  equality defines $\lambda$.
Then $  \eta \sim \ln \lambda^{-1} $, and
\begin{align*}
\mu(\Omega_\eta) \sim e^{2\eta}\, \eta^{\frac{n-1}2} =  \lambda^{-1}\,
\eta^{\frac{q-n-1}2} \eta^{\frac{n-1}2}
\sim  \lambda^{-1}\, (\ln \lambda^{-1} )^{\frac{q}2-1}.
\end{align*}
For $q > 2$  this shows that   \eqref{rf1}   is sharp, and in particular
that $R_D$  is not of weak type
$(1,1)$.

Theorem   \ref{one'} is completely proved.
\hfill $\Box$

\renewcommand{\theequation}{\thesection.\arabic{equation}}
\section{Proof of Theorem \ref{H} (a) } \label{sec4}
\setcounter{equation}{0}

\medskip

Before proving this theorem, we find local and global estimates for
$t^{{k}/{2}} D_x p_t(x, y)$
in the space $L^2(dt/t)$. They will be analogous to Propositions \ref{local1}  and   \ref{glob2},
with similar proofs.

 \begin{propo}\label{local5}
For $0 < | x - y | \le 1$, one has
\begin{equation}
  \label{stest51}
\|t^{\frac{k}{2}} D_x\, p_t(x, y)\|_{L^2(dt/t)}  \lesssim \frac{1}{\mu(B(x, | x - y |))}
\end{equation}
and
\begin{equation}  \label{stest52}
 \|t^{\frac{k}{2}} \nabla_{y}\,D_x \,p_t(x, y)\|_{L^2(dt/t)}  \lesssim
\frac{1}{|x-y|\,\mu(B(x, | x - y |))}.
\end{equation}
 \end{propo}

\begin{propo} \label{glob52}
For $ | x - y | > 1$, one has
 \begin{equation} \label{glob51}
  \|t^{\frac{k}{2}} D_x\, p_t(x, y)\|_{L^2(dt/t)}  \lesssim  
e^{-x_1-y_1 -|x-y|}\, |x-y|^{\frac{q-n}2-\frac14}
\left[1 + \left(\frac{|x'-y'|^2}{|x-y|}  \right)^{\frac{k}2}\right].
 \end{equation}
\end{propo}

\vskip12pt

\noindent {\bf Proof of Propositions \ref{local5} and \ref{glob52}.}
Clearly,
\begin{align*}
\left|t^{\frac{k}{2}} D_x\, p_t(x, y)\right|^2 = \sum_{|\alpha| = |\beta|=k}
a_\alpha\, \overline{a_\beta}\, t^k \partial_x^\alpha p_t(x, y)\, \partial_x^\beta p_t(x, y).
\end{align*}
We write  $\alpha = (\alpha_1, \alpha')$ as before and similarly
 $\beta = (\beta_1, \beta')$, assuming always $|\alpha| = |\beta| = k$.
Using the expression \eqref{KLF1} for $p_t(x, y)$ and  \eqref{deriv},
we see that  $t^k \partial_x^\alpha \,p_t(x, y)\, \partial_x^\beta \,p_t(x, y)$
 is the sum of positive factors times
\[
 e^{-2x_1-2y_1} \,
t^{k-\frac{|\tilde \alpha|+|\tilde \beta|}2-n}\,
H_{\tilde \alpha}\left(\frac{x - y}{2 \sqrt{t}}\right)\,
H_{\tilde \beta}\left(\frac{x - y}{2 \sqrt{t}}\right)\,
 e^{- 2t}\,e^{-\frac{|x-y|^2}{2 t}},
\]
taken over  $ \tilde \alpha_1 = 0, \dots, \alpha_1$ and
 $ \tilde \beta_1 = 0, \dots, \beta_1$.
Here as in  Section \ref{sec2}, we have
$\tilde \alpha = ( \tilde \alpha_1,\alpha') $, and similarly
$\tilde \beta = (\tilde \beta_1, \beta')$.

Integrating this expression with respect to  $dt/t$,  
we obtain a positive constant times
\begin{align}\label{term2}
e^{-2x_1-2y_1} \,
\int_0^{+\infty} t^{k-\frac{|\tilde \alpha|+|\tilde \beta|}2-n}\,
H_{\tilde \alpha}\left(\frac{x - y}{2 \sqrt{t}}\right)\,
H_{\tilde \beta}\left(\frac{x - y}{2 \sqrt{t}}\right)\,
 e^{- 2t}\,e^{-\frac{|x-y|^2}{2 t}} \,\frac{dt}t.
\end{align}
Here the product of the two Hermite polynomials is a linear combination
of terms
\begin{align*}
\left(\frac{x_1 - y_1}{2 \sqrt{t}}\right)^{ \gamma_1}
\left(\frac{x' - y'}{2 \sqrt{t}}\right)^{\gamma'}
\left(\frac{x_1 - y_1}{2 \sqrt{t}}\right)^{\delta_1}
\left(\frac{x' - y'}{2 \sqrt{t}}\right)^{\delta'},
\end{align*}
where  $\gamma = (\gamma_1, \gamma')$ and   $\delta =  (\delta_1, \delta')$
satisfy $0 \le \gamma_1 \le \tilde \alpha_1$ and
$0 \le \gamma' \le \alpha'$,  and similarly
$0 \le \delta_1 \le \tilde \beta_1$ and  $0 \le \delta' \le  \beta'$.
Each such term leads to a
 corresponding term in the expression   \eqref{term2},  proportional to
\begin{align}\label{term3}
e^{-2x_1-2y_1} \, (x_1 - y_1)^{\gamma_1+\delta_1} (x' - y')^{\gamma'+\delta'}
\int_0^{+\infty} t^{k-\frac{|\tilde \alpha|+|\tilde \beta|+|\gamma|+|\delta|}2-n}\,
 e^{- 2t}\,e^{-\frac{|x-y|^2}{2 t}} \,\frac{dt}t.
\end{align}
After replacing $t$ by $\tau = 2t$ in this integral, one can apply Lemma
\ref{intest} with $a = 2|x-y|$.

For $|x - y| \le 1$, we see from  Lemma \ref{intest} (a)  that the modulus
of  \eqref{term3} is controlled by
\begin{align*}  
e^{-2x_1-2y_1} \, |x_1 - y_1|^{\gamma_1+\delta_1} |x' - y'|^{|\gamma'|+|\delta'|}
|x-y|^{2k-(|\tilde \alpha|+|\tilde \beta|+|\gamma|+|\delta|)-2n} \lesssim
e^{-4x_1}\,|x-y|^{-2n},
\end{align*}
if  the exponent of $t$ in the
integral is negative. But in the opposite case, this is also true.
 Thus  \eqref{stest51}
is proved. To complete the proof of Proposition \ref{local5}, one verifies
\eqref{stest52} in a similar way.

Aiming at Proposition \ref{glob52}, we assume that $|x - y| > 1$ and
conclude from Lemma \ref{intest}~(b)
that  \eqref{term3} is in modulus bounded by constant times
\begin{align*} 
   e^{-2x_1-2y_1-2|x-y|} \, |x_1 - y_1|^{\gamma_1+\delta_1} |x' - y'|^{|\gamma'|+|\delta'|}
|x-y|^{k-\frac{|\tilde \alpha|+|\tilde \beta|+|\gamma|+|\delta|}2-n-\frac12}.
\end{align*}
In this expression, we estimate $|x_1 - y_1|$ by  $|x - y|$ and write $k$ as
$k =  (\alpha_1 + |\alpha'| +\beta_1 + |\beta'|)/2$,
getting
\begin{align} \label{term4}
   e^{-2x_1-2y_1-2|x-y|} \,
|x-y|^{\frac{\alpha_1 +\beta_1}2-\frac{\tilde \alpha_1 -\gamma_1}2-\frac{\tilde \beta_1 -\delta_1}2-n-\frac12}\left(\frac{|x'-y'|}{|x-y|^{1/2}}  \right)^{|\gamma'|+|\delta'|}.
\end{align}
Here we have
\begin{equation}
  \label{inequ}
  \alpha_1 \le q, \quad \beta_1 \le q, \quad
\tilde \alpha_1 - \gamma_1\ge 0, \quad   \tilde \beta_1 - \delta_1 \ge 0
\end{equation}
 and also $|\gamma'|,\;|\delta'| \le k$.
 Thus
for $|x - y| \ge 1$
the last expression is at most
\begin{align} \label{term8}
   e^{-2x_1-2y_1-2|x-y|} \, |x-y|^{q-n-\frac12}
\left[1 + \left(\frac{|x'-y'|^2}{|x-y|}  \right)^{k}\right].
\end{align}
We have shown that this is a bound for $\int_0^{+\infty} t^k |D_x \,p_t(x,y)|^2\,dt/t$,
and \eqref{glob51} follows.

The two propositions are proved.  \hfill $\Box$

\vskip7pt

\noindent \textbf{Remark} The estimate at the end of the proof just given can be
sharp only when all the inequalities \eqref{inequ} are equalities. Indeed,
if any of the four inequalities is strict, one can estimate \eqref{term4}
by a modified quantity   \eqref{term8} in which $|x-y|$ has a smaller
exponent.

\vskip7pt

\noindent \textbf{Proof of Theorem \ref{H} (a).}
Clearly,     
\begin{align*}
\mathcal{H}_{D}(f)(x) = \left( \int_0^{+\infty} \left| t^{\frac{k}{2}} \int D_x\, p_t(x, y) f(y) \, d\mu(y) \right|^2 \, \frac{dt}{t} \right)^{\frac{1}{2}}.
\end{align*}
We now introduce a local part $\mathcal{H}_{D}^{\mathrm{loc}}$ and a global
part  $\mathcal{H}_{D}^{\mathrm{\infty}}$ of this operator, defined by
restricting the integration in $y$ here to the region $|y-x| \le 1$ and
$|y-x| > 1$, respectively. Thus
\begin{align*}
\mathcal{H}_{D}f \le \mathcal{H}_{D}^{\mathrm{loc}}f +
 \mathcal{H}_{D}^{\infty}f.
\end{align*}

To deal with the local part, we apply the method of localization as
done for the Riesz transforms, but now using vector-valued
singular integral theory and  Proposition \ref{local5}.

For  $\mathcal{H}_{D}^{\infty}$, we start by applying Minkowski's integral
inequality, getting
\begin{align*}
\mathcal{H}_{D}^{\infty}(f)(x) \le  \int_{|x - y| > 1} |f(y)| \,
\left( \int_0^{+\infty} \left| t^{\frac{k}{2}} D_x \,p_t(x, y) \right|^2 \frac{dt}{t}  \, \right)^{\frac{1}{2}} \, d\mu(y).
\end{align*}

We can now follow the argument for  $R_D^{\infty}$ in the proof of Theorem \ref{one'}, replacing  \eqref{glob1} by \eqref{glob51}. In particular, that part of $\mathcal{H}_{D}^\infty$
given by  $x_1-y_1 \le 1$ is seen to be of strong type $(1,1)$.
Moreover, the kernel of the
remaining part is controlled by $\mathcal V_{\kappa}(x,y)$
with $\kappa = q+1/2.$ This leads to the weak type $(1,1)$ estimate
of  Theorem  \ref{H} (a) and also to  \eqref{rf5'},
in view of Propositions \ref{prop2} and \ref{propn1}.


For the sharpness parts of Theorem  \ref{H} (a), we let $f = \chi_{B(0,1)}$.
Let   $|\alpha| = k$ and take points $y$ and $x$ with
 $|y|<1$ and  $x_1>0$  large but $x'/\sqrt{|x|}$  bounded, which implies
$x_1 = |x|+ O(1)$. We will then estimate
 $t^{k/2}\partial_x^\alpha\, p_t(x,y)$, and we need to be more precise than in
the preceding argument.
Because of \eqref{KLF1}  and \eqref{deriv},
this quantity is the sum of   
\begin{align}   \label{termP}
(-1)^{k}\,\binom{\alpha_1}{\tilde \alpha_1}
\,2^{-|\tilde \alpha|}e^{-x_1-y_1}\,t^{\frac{k-|\tilde \alpha|-n}2}\, H_{\tilde \alpha}\left(\frac{x-y}{2\sqrt t}\right)
e^{-t} e^{-\frac{|x-y|^2}{4 t}},
\end{align}
  taken over
$\tilde \alpha_1 = 0,\dots, \alpha_1$. Here
$\tilde \alpha = ( \tilde \alpha_1,\alpha') $
as before. The  remark after the end of the proof of Propositions
\ref{local5} and  \ref{glob52} shows that we need only consider the case
$ \alpha_1 = q$, which implies
$| \alpha'| = k-q$.
For the same reason, we may replace the polynomial $H_{\tilde \alpha_1}\left({(x_1-y_1)}/{2\sqrt t}\right) $
by its leading term
$2^{\tilde \alpha_1}\left({(x_1-y_1)}/{ 2 \sqrt{t}}\right)^{\tilde \alpha_1}$.

Then we can write
\begin{align*}   
 H_{\tilde \alpha}\left(\frac{x-y}{2\sqrt t}\right)
=&\, 2^{ \tilde \alpha_1}  \left(\frac{x_1-y_1}{2\sqrt t}\right)^{\tilde \alpha_1}
H_{ \alpha'}\left(\frac{x'-y'}{2\sqrt t}\right)+\dots \\
=&\, 2^{ \tilde \alpha_1}\,  t^{\frac{{\tilde \alpha_1}}{2}}
 \left(\frac{x_1-y_1}{2t}\right)^{\tilde \alpha_1}
H_{ \alpha'}\left(\frac{x'-y'}{2\sqrt t}\right)+\dots,
\end{align*}
where the dots indicate negligible terms.
Restricting  $t$  by $\left|t - \frac{|x|}2\right| < {|x|^{3/4}}$, we
see that
\begin{align*}
 H_{\tilde \alpha}\left(\frac{x-y}{2\sqrt t}\right)
=&\, 2^{ \tilde \alpha_1}\,  t^{\frac{{\tilde \alpha_1}}{2}}
\left(\frac{x_1}{|x|} + O\left(|x|^{-\frac14}\right) \right)^{\tilde \alpha_1}
H_{ \alpha'}\left(\frac{x'}{\sqrt{2|x|} }+ O\left(|x|^{-\frac14}\right) \right)
+\dots\\ =&
\, 2^{ \tilde \alpha_1} \, t^{\frac{{\tilde \alpha_1}}{2}}
H_{ \alpha'}\left(\frac{x'}{\sqrt{2|x|} }\right)
+ \dots.    
\end{align*}

We now insert this last expression for the Hermite polynomial in \eqref{termP},
and observe that $|\tilde \alpha| =\tilde \alpha_1 +| \alpha'|= \tilde \alpha_1 +k-q $. Then the only
factor in \eqref{termP} which depends on $\tilde \alpha_1$ will be the
binomial coefficient. Since
\begin{align*}
\sum_{\tilde \alpha_1 = 0}^{\alpha_1} \,\binom{\alpha_1}{\tilde \alpha_1}  =
2^{\alpha_1} = 2^q,
\end{align*}
we get by summing  \eqref{termP} in $\tilde \alpha_1$
\begin{align*}
t^{\frac k2}\,\partial_x^\alpha\, p_t(x,y) =
(-1)^{k}\,
\,2^{2 q - k} e^{-x_1-y_1}\,t^{\frac{q-n}2}\, H_{ \alpha'}\left(\frac{x'}{\sqrt{2|x|} }\right)
e^{-t} e^{-\frac{|x-y|^2}{4 t}} + \dots.
\end{align*}

The next step is to sum these expressions in $\alpha$,
with the coefficients from
 \eqref{D}.  The result is that for $y,\:x$ and $t$  as described above,
\begin{align*}
 &t^{\frac k2}\,D_x\, p_t(x,y)\\& =  \,  (-1)^{k}\,2^{2 q - k}\,e^{-x_1-y_1}\,
\, t^{\frac{q - n}2} \,  e^{-t}\,  e^{-\frac{|x-y|^2}{4 t}}
\sum_{|\alpha'|=k-q} a_{(q,\alpha')}\,H_{ \alpha'}\left(\frac{x'}{\sqrt{2|x|} }\right)
+\mathrm{harmless\hskip5pt terms}.
\end{align*}

The sum here is $P\left(x'/\sqrt{2|x|}\right)$ for some  nonzero polynomial $P$ in
$n-1$ variables. Thus we can find a ball $B\subset \R^{n-1}$ in which  $P$
does not vanish.
Further,
\begin{align*}   
e^{-\frac{|x-y|^2}{4 t}} \sim e^{-\frac{|x|^2}{4 t}},
\end{align*}
since  $|y|<1$ and $t \sim |x|$. With $f = \chi_{B(0,1)}$, we then get by
integrating in $y$
\begin{align*}
\left| t^{\frac{k}{2}} D\,e^{t \Delta_{e_1}} f(x) \right|^2 \sim
e^{-2x_1}\,x_1^{q-n}  e^{-2t}e^{-\frac{|x|^2}{2 t}}
\end{align*}
 where as before  $x_1>0$  is large, $x'/\sqrt{2 { |x|}} \in B$ and
$\left|t - \frac{|x|}2\right| < {|x|^{3/4}}$. Hence,
\begin{align}  \label{HD}
 \int_0^{+\infty}
\left| t^{\frac{k}{2}} D\,
e^{t \Delta_{e_1}} f(x) \right|^2 \frac{dt}{t} \gtrsim
e^{-2x_1}\,x_1^{q-n}\, \int_{\frac{|x|}2-|x|^{3/4}}^{\frac{|x|}2+|x|^{3/4}}
\, e^{-2t}e^{-\frac{|x|^2}{2 t}}
 \frac{dt}{t}.
\end{align}
After the change of variable $\tau = 2t$, we can apply the last part of
 Lemma  \ref{intest} (b)
to conclude that the right-hand side of \eqref{HD} equals a positive constant times
\begin{align}\label{HH}
e^{-2x_1}\,x_1^{q-n}\, e^{-2|x|}\,|x|^{-\frac{1}{2}}\, (1+O(|x|^{-1}))
= e^{-2x_1-2|x|}\,x_1^{q-n-\frac{1}{2}}\, (1+O(|x|^{-1})).
\end{align}

We define  $\Omega_\eta$ by \eqref{omeg}, though $B$ is
 not the same as in  Section~\ref{sec3}. If $\eta$ is large,
$x \in \Omega_\eta$ implies  $x'/\sqrt{2 { |x|}} \in B$.
It now follows from   \eqref{HD} and  \eqref{HH} that
for large $\eta$
\begin{align*}    
\mathcal{H}_{D}^\infty f
\gtrsim  e^{-2\eta}\, \eta^{\frac {q-n}2-\frac{1}{4}}.
\end{align*}
 in $\Omega_\eta$.
This implies the ``only if'' part of Theorem  \ref{H} (a)
and the sharpness of $q$;
 cf.\ the very last part of Section  \ref{sec3}. We leave it
to the reader to check that all the neglected terms above can be disregarded.

 Theorem  \ref{H} (a) is  proved.
 \hfill $\Box$

\renewcommand{\theequation}{\thesection.\arabic{equation}}
\section{Proof of Theorem \ref{H} (b)} \label{sec5}
\setcounter{equation}{0}

\medskip

To prove the claimed boundedness properties of  $\mathcal{G}_{\alpha}$,
we copy the reasoning used for Theorems  \ref{one'} and \ref{H} (a). We
must  estimate the kernel $t^k\,D_xP_t(x,y)$ in the space $L^2(dt/t)$.

The subordination formula says that the Poisson kernel $P_t$ is given by
\begin{align*}
P_t(x, y) = \frac{t}{2 \sqrt{\pi}}
\int_0^{+\infty} u^{-\frac{3}{2}}\, e^{-\frac{t^2}{4 u}}\, p_u(x, y) \, du.
\end{align*}
This will also hold with  $P_t(x, y)$ and $p_u(x, y)$ replaced by their
derivatives  $D_xP_t(x, y)$ and $D_xp_u(x, y)$, respectively. Using
Minkowski's integral inequality, we conclude that
\begin{align*}
\left(\int_0^{+\infty} |t^kD_xP_t(x,y)|^2\,\frac{dt}t\right)^{\frac12}
&\sim
\left(\int_0^{+\infty} \left|t^{k+1}\int_0^{+\infty} u^{-\frac{3}{2}}\, e^{-\frac{t^2}{4 u}}\, D_xp_u(x, y) \, du\right|^2\,\frac{dt}t\right)^{\frac12}\\
&\le \int_0^{+\infty}  u^{-\frac{3}{2}} |D_xp_u(x, y)| \left(\int_0^{+\infty}
\left(t^{k+1}\, e^{-\frac{t^2}{4 u}}\right)^2 \,\frac{dt}t\right)^{\frac12} \, du.
\end{align*}
The last inner integral here is
   seen to equal $2^{k}\,\Gamma(k+1) u^{k+1}$,
and so
\begin{align} \label{DxP}
\left(\int_0^{+\infty} |t^k\,D_xP_t(x,y)|^2\,\frac{dt}t\right)^{\frac12}
\ls\int_0^{+\infty}  u^{\frac{k}{2}} |D_xp_u(x, y)|\,\frac{du}u.
\end{align}

The right-hand side in \eqref{DxP} is like the expression for
$R_D(x,y)$ in \eqref{RR}, except for the modulus signs. In Section \ref{sec2},
the equality \eqref{RR} was used to prove Propositions \ref{local1}
and   \ref{glob2}.
 Tracing those arguments, one sees that
the inequality  \eqref{DxP} is sufficient to imply estimates similar to
\eqref{stest1} and \eqref{glob1} for
$(\int_0^{+\infty} |t^k\,D_xP_t(x,y)|^2\,{dt}/t)^{1/2}$.
Taking a derivative also in $y$, one  obtains the analog of \eqref{stest2}.
This leads to  the  weak type of   $\mathcal{G}_{\alpha}$ for $q\le 2$
and also to  \eqref{rf4}, as seen from the proof of Theorem $\ref{one'}$
in Section $\ref{sec3}$.

It remains to verify the sharpness parts of Theorem \ref{H} (b).
Recall that $D = \sum_{|\alpha| = k} a_{\alpha} \partial^{\alpha}$.
We will find an
 $f \in L^1(\mu)$
and a  ball $B \subset \rr^{n-1}$ such that if $x_1 > 0$ is large and $x'/\sqrt{x_1} \in B$, one has
 \begin{align}
 \mathcal{G}_Df(x) \gtrsim e^{-2x_1}\,x_1^{\frac{q-n-1}2}. \label{lowerest}
\end{align}
From this estimate, the sharpness parts of Theorem  \ref{H} (b) will follow,
cf.\ the last few lines of Section \ref{sec3}.

To estimate $D_x P_t(x, y)$, we consider
\begin{align*}
\partial^\alpha_x P_t(x, y) = \frac{t}{2 \sqrt{\pi}}
\int_0^{+\infty} u^{-\frac{3}{2}}\, e^{-\frac{t^2}{4 u}}\,
\partial^\alpha_x  p_u(x, y) \, du,
\end{align*}
with $|\alpha| = k$ and  $\alpha_1 \le q$.

    From \eqref{KLF1}
we see   that $\partial^\alpha_x P_t(x, y)$
is a sum of terms   
\begin{align*}   
(-1)^{\alpha_1-\tilde \alpha_1}\, (4\pi)^{-\frac {n+1} 2}\, \binom{\alpha_1}{\tilde \alpha_1} \,e^{-x_1-y_1}\, t
\int_0^{+\infty} u^{-\frac{1+n}{2}}\, e^{-\frac{t^2}{4 u}}\, e^{-u}\,\partial^{\tilde \alpha}_x  \, e^{-\frac{|x-y|^2}{4 u}}\,
\frac{du}u,
\end{align*}
where
 $\tilde \alpha = (\tilde \alpha_1, \alpha')$ as before,
and the sum is taken over $ \tilde  \alpha_1 = 0,\dots,\alpha_1$.
Using \eqref{deriv},  we see that this expression equals
 \begin{multline} \label{term88}
(-1)^{k}\, \pi^{-\frac {n+1} 2}\,2^{-|\tilde \alpha|-n-1}\,\binom{\alpha_1}{\tilde \alpha_1}\, e^{-x_1-y_1}\, t
 \int_0^{+\infty}
u^{-\frac{1 +|\tilde \alpha|+n}2}  H_{\tilde  \alpha}\left(\frac{x - y}{2 \sqrt{u}}\right) e^{- u}\,
e^{-\frac{t^2+|x-y|^2}{4 u}} \,\frac{du}u.
\end{multline}

We now apply \eqref{ext} to the integral here.
It follows that for large $|x-y|$ the quantity  \eqref{term88} equals
 \begin{multline}   \label{term888}
(-1)^{k}\, (2\pi)^{-\frac {n} 2}\,2^{-\frac{|\tilde \alpha|} 2}\,
\binom{\alpha_1}{\tilde \alpha_1}\,  t\,
(t^2+|x-y|^2)^{-\frac{1 +|\tilde \alpha|+n}4-\frac{1}4}\, e^{-x_1-y_1-\sqrt{t^2+|x-y|^2}}\\ \times
  H_{\tilde  \alpha}\left(\frac{x - y}{\sqrt 2\, (t^2+|x-y|^2)^{1/4}}\right)
\left(1+O\left(\frac1{\sqrt{t^2+|x-y|^2}}\right)\right).
\end{multline}

Again we
 consider points $x$
  for which $x_1 > 0$ is large and $x'/\sqrt{x_1}$   stays bounded.
Further, we assume $|y|<1$ and restrict $t$ by $\sqrt{x_1}< t < 2\sqrt{x_1}$.
Some Taylor expansions then lead to
$|x-y| = x_1 \left(1+O({x_1}^{-1})\right)$ and
\begin{align}
  \sqrt{t^2+|x-y|^2}\, =&\, |x|\left(1 - \frac{x\cdot y}{|x|^2} + \frac{t^2}{2|x|^2}
+ O\left(\frac1{|x|^2}\right) \right)  \notag \\ =& \,
 |x| - y_1 + \frac{t^2}{2x_1} +  O\left(\frac1{\sqrt{x_1}}\right)
\label{taylor}\\  =& \,x_1 \left( 1 + O(x_1^{-1}) \right), \notag
\end{align}
as $ x_1 \to \infty$.

Using \eqref{taylor}, we conclude that the expression \eqref{term888} now equals
 \begin{multline}   \label{term21}
(-1)^{k}\, (2\pi)^{-\frac {n} 2}\,2^{-\frac{|\tilde \alpha|} 2}\,
\binom{\alpha_1}{\tilde \alpha_1}\, t \,
x_1^{-\frac{|\tilde \alpha|+n+2}2}\, e^{-x_1 - |x|  - \frac{t^2}{2x_1}}\,\\ \times
  H_{\tilde  \alpha_1}\left(\sqrt{\frac{x_1}{2}}\right)
H_{\alpha'}\left(\frac{x'}{\sqrt{2x_1}}\right)
 \left(1+ O\left(\frac1{\sqrt{x_1}}\right)\right),
\end{multline}
as $ x_1 \to \infty$. We can replace the Hermite polynomial
$H_{\tilde  \alpha_1}\left(\sqrt{{x_1}/{2}}\right)$ here by its leading term,
which is $2^{\tilde  \alpha_1/2}x_1^{\tilde  \alpha_1/2}$. Since $|\tilde \alpha| =
\tilde  \alpha_1 + | \alpha'| $, this will make
the exponents of $2$ and  $x_1$ in  \eqref{term21} independent of
$\tilde  \alpha_1$. When we then sum in $\tilde  \alpha_1 = 0,\dots, \alpha_1$,
the binomial coefficients in \eqref{term21} will sum up to $2^{ \alpha_1}$,
and as a result
\begin{equation*}
  \partial^\alpha_x P_t(x, y) = (-1)^{k}\, (2\pi)^{-\frac {n} 2}\,
2^{ \alpha_1-\frac{|\alpha'|}2}\, t \,
x_1^{-\frac{|\alpha'|+n+2}2} \, e^{-x_1 - |x|  - \frac{t^2}{2x_1}}\,
H_{\alpha'}\left(\frac{x'}{\sqrt{2x_1}}\right)
\left(1+ O\left(\frac1{\sqrt{x_1}}\right)\right).
\end{equation*}
Next, we integrate against $f(y)\,d\mu(y)$ with $f = \chi_{B(0,1)}/\mu(B(0,1))$.
Summing also over $\alpha$, we get
\begin{multline*}
  D_x P_tf(x) \\
= (-1)^{k}\, (2\pi)^{-\frac {n} 2}\,
 \, t \,x_1^{-\frac{n+2}2} \,e^{-x_1 - |x|  -\frac{t^2}{2x_1}}\,
\sum_{\alpha}2^{ \alpha_1-\frac{|\alpha'|}2}  a_{\alpha}\,
x_1^{-\frac{|\alpha'|}2} H_{\alpha'}\left(\frac{x'}{\sqrt{2x_1}}\right)
\left(1+ O\left(\frac1{\sqrt{x_1}}\right)\right).
\end{multline*}

For large $x_1$, the largest  terms in the sum here are those where $|\alpha'|$ takes its
minimal value $k-q$, so that $\alpha_1 = q$.  We introduce the nonzero polynomial
\begin{equation*}
 P(z') =  (-1)^{k}\, (2\pi)^{-\frac {n} 2}\, 2^{\frac{3q-k}2}
\sum_{|\alpha'|=k-q} a_{(q,\alpha')}
 H_{\alpha'}(z').
\end{equation*}

Since $x'/\sqrt{2x_1}$ stays bounded, it follows that
\begin{equation*}
 t^k\, D_x P_tf(x) =
 e^{-x_1 - |x|  - \frac{t^2}{2x_1}}\, t^{k+1} \,x_1^{-\frac{k-q+n+2}2}
 \left( P\left(\frac{x'}{\sqrt{2x_1}}   \right)+
 O\left(\frac1{\sqrt{x_1}}\right)  \right).
\end{equation*}

We can now integrate in $t$ and conclude that
\begin{equation*}
 \left(\mathcal G_D f(x)\right)^2  \gtrsim   e^{-4x_1 }
\int_{\sqrt{x_1}}^{2\sqrt{x_1}}
 e^{ -\frac{t^2}{x_1}}\, t^{2k+2}\,\frac{dt}t \;\;x_1^{-({k-q+n+2})}
 \left( P\left(\frac{x'}{\sqrt{2x_1}}   \right)^2 +
 O\left(\frac1{\sqrt{x_1}}\right)  \right).
\end{equation*}
The integral here is of order of magnitude $x_1^{k+1}$, and so
\begin{equation*}
 \mathcal G_D f(x) \gtrsim   e^{-2x_1 }x_1^{\frac{q-n-1}2}
 \left( \left|P\left(\frac{x'}{\sqrt{2x_1}}\right)\right|   +
 O\left(\frac1{\sqrt{x_1}}\right)  \right).
\end{equation*}

 To obtain    \eqref{lowerest},   it is now
enough to let  $B$ be a closed ball in    $\rr^{n-1}$  in
which the polynomial $P$ does not vanish.

This ends the proof of Theorem \ref{H} (b).
\hfill $\Box$

\medskip

\renewcommand{\theequation}{\thesection.\arabic{equation}}
\section{Sharp estimates for $\frac{\partial^k}{\partial t^k} \,p_t$}\label{sec6}
\setcounter{equation}{0}

\medskip

This section is a preparation for the proof of Theorem \ref{TH2}.

We take $k$ derivatives with respect to $t$ of
\begin{align} \label{pt}
p_t(x, y) = (4 \pi)^{-\frac{n}{2}} e^{-x_1 - y_1} t^{-\frac{n}{2}} e^{-t - \frac{|x - y|^2}{4 t}}.
\end{align}
Those derivatives which fall on the last exponential here will produce a factor $\frac{|x - y|^2}{4 t^2} - 1$. The derivative
\begin{align*}
\frac{\partial}{\partial t} \left[  \frac{|x - y|^2}{4 t^2} - 1 \right] = -\frac{|x - y|^2}{2 t^3}
\end{align*}
will also appear.
We see that
\begin{align}\label{prod}
\frac{\partial^k}{\partial t^k}\, p_t(x, y) = q_k\, p_t(x, y),
\end{align}
where the factor $q_k$ is given by
\begin{align*}
q_k = Q_k\left(\frac{|x - y|^2}{4 t^2} - 1, \frac{1}{t}, \frac{|x - y|^2}{2 t^3}\right)
\end{align*}
for a polynomial $Q_k$ in three variables, whose coefficients depend only on $n$ and $k$.

To estimate $q_k$, we examine the terms of this polynomial. Consider a term obtained by letting exactly $m$ of the $k$ differentiations fall on the exponential $e^{-t - {|x - y|^2}/{4 t}}$. If, moreover, the number of  differentiations falling on a power of $\frac{|x - y|^2}{4 t^2} - 1$ is $p$, the resulting term of $Q_k$ will contain $p$ factors $-\frac{|x - y|^2}{2 t^3}$ and $m - p$ factors $\frac{|x - y|^2}{4 t^2} - 1$. The remaining $k - m -  p$ differentiations will produce factors $-t^{-1}$. We conclude that  $q_k$ is a sum of expressions of the form
\begin{align} \label{6Y}
C_{m, p}  \left(\frac{|x - y|^2}{4 t^2} - 1\right)^{m - p}  \left( -\frac{|x - y|^2}{2 t^3} \right)^p (-t)^{-(k - m - p)}
\end{align}
with $C_{m, p}  > 0$. Here $p + m \leq k$ and $p \leq m$; thus $2 p \leq k$ so that $p \leq [\frac{k}{2}]$. This implies the upper estimate
\begin{align} \label{6A}
\left| \frac{\partial^k}{\partial t^k}\, p_t(x, y) \right| \lesssim \sum \left| \frac{|x - y|^2}{4 t^2} - 1 \right|^{m - p} |x - y|^{2 p}\, t^{-k + m - 2p}\,  p_t(x, y),
\end{align}
 the sum taken over  $0 \le p \le [\frac k 2]$ and  $0 \le m \le k-p$.

We also need a lower estimate.
Let $\eta > 0$ be large and define
\begin{align*}
\Sigma_{\eta} = \{x;\; \eta - 1 < x_1 < \eta \; \; \;  \mbox{and } \;  \; \sqrt{\eta} < x_i < 2 \sqrt{\eta},\;\; i = 2, \cdots, n \}.
\end{align*}

\begin{lem} \label{lemzz}
For $x \in \Sigma_{\eta}$,\hskip8pt $|y| < 1$ and
\begin{align*}
\frac{\eta}{2} \left(1 - 2 \frac{c_1}{\sqrt{\eta}}\right) < t  < \frac{\eta}{2} \left(1 -  \frac{c_1}{\sqrt{\eta}}\right)
\end{align*}
with $c_1 = c_1(n, k)$ small enough and $\eta$ large enough, one has
\begin{align}  \label{zzs6}
(-1)^{[\frac{k}{2}]}\, t^k \frac{\partial^k}{\partial t^k}\, p_t(x, y)  \gtrsim e^{-2 \eta}\, \eta^{\frac{k - n}{2}}.
\end{align}
\end{lem}

\noindent \textbf{Proof.}  Assuming $x,\; y$ and $t$ as in the lemma, we can write
\begin{equation*}
  t =  \frac{\eta}{2}\left(1 -  \frac{\zeta}{\sqrt{\eta}}\right),
\end{equation*}
where  $c_1 < \zeta < 2 c_1$.

Some simple computations will lead to
$|x - y| = \eta \left( 1 + O(\eta^{-1}) \right)$
and
\begin{align*}
\frac{|x - y|^2}{4 t^2} - 1 = 2\, \frac{\zeta}{\sqrt{\eta}}  \left( 1 + O(\eta^{-\frac{1}{2}}) \right) \qquad  \mbox{and} \qquad
\frac{|x - y|^2}{2 t^3} = \frac{4}{\eta}  \left( 1 + O(\eta^{-\frac{1}{2}}) \right),
\end{align*}
 as $\eta \to \infty$.
The expression \eqref{6Y} thus equals
\begin{multline}\label{Q}
(-1)^{m - k}\, C_{m, p} \left( 2 \frac{\zeta}{\sqrt{\eta}} \right)^{m - p} \left( \frac{4}{\eta} \right)^p \left( \frac{\eta}{2} \right)^{-k + m + p} \left(1 + O(\eta^{-\frac{1}{2}} \right)\\ = (-1)^{m - k}\, C_{m, p}\, 2^k\, \eta^{\frac{m + p}{2} - k}\, \zeta^{m - p} \left(1 + O(\eta^{-\frac{1}{2}}) \right).
\end{multline}

If we fix a small $c_1>0$ and take $\eta > c_1^{-1}$, the product $ 2^k\, \eta^{\frac{m + p}{2} - k}\, \zeta^{m - p}$ will be maximal  when $m + p$ takes its largest possible
value $k$ and also $m - p$ is as small as possible. This means that $p = [\frac{k}{2}]$ and  $m = k - [\frac{k}{2}]$. Moreover, if  $c_1$ is small enough and  $\eta$
large enough, these values of $m$ and $p$ will make the expression \eqref{Q} much
larger than any other admissible values of  $m$ and $p$, in absolute value.
We conclude that then
\begin{align}\label{q}
 (-1)^{[\frac{k}{2}]}\, q_k \sim \eta^{-\frac{k}{2}} \zeta^{k - 2 [\frac{k}{2}]} \sim \left\{\begin{array}{rl}
\eta^{-\frac{k}{2}}, & k \mbox{ even} \\
c_1 \eta^{-\frac{k}{2}}, & k \mbox{ odd.}
\end{array}
\right.
\end{align}

We also need an estimate of the value of $p_t(x,y)$ for these $x$,  $y$,  $t$,
and write
the sum of the exponents of $e$ in \eqref{pt} as
\begin{align*}
  -x_1 - y_1 -t - \frac{|x - y|^2}{4 t} =& \,-\eta +O(1) - \frac\eta2
 +\frac{\zeta\sqrt\eta}2 - \frac{\eta^2(1+O(\eta^{-1}))}{2\eta(1-\zeta/\sqrt \eta)}\\
=&\, -\eta +O(1) - \frac\eta2 +\frac{\zeta\sqrt\eta}2 - \frac\eta2\left(1+\frac\zeta{\sqrt \eta}
+O(\eta^{-1})\right)  \\
=& -2\eta +O(1).
\end{align*}
Thus  \eqref{pt} implies that
\begin{align*}
 p_t(x,y) \sim e^{-2\eta}\, t^{-\frac n2}.
\end{align*}
We now combine this with \eqref{prod} and \eqref{q}, and take
 $\eta$ large after fixing a small $c_1$. Then \eqref{zzs6} follows.
\hfill $\Box$

\medskip

\renewcommand{\theequation}{\thesection.\arabic{equation}}
\section{Proof of Theorem \ref{TH2}} \label{sec7}
\setcounter{equation}{0}

\medskip

For the local parts of the operators in Theorem  \ref{TH2},
one can use the method of
localization, since $\R^{(n, v)}$ has the local doubling property.
 Standard vector-valued singular integral  theory then
 gives the weak type $(1, 1)$ of the local parts of $g_k$, $h_k$ and $H_k$.

The parts at infinity of these operators can all be estimated by the
method used for Theorem \ref{one'}, in the following way.

Consider first the part at infinity of $h_k$, given by
\begin{align*}
h_k^{\infty}(f)(x) = \left( \int_0^{+\infty}  \left|t^{k} \int_{|x - y|
> 1}  \frac{\partial^k}{\partial t^k}\, p_t(x, y)\, f(y) \,
d\mu(y) \right|^2 \, \frac{dt}t \right)^{\frac{1}{2}}.
\end{align*}
Using Minkowski's integral inequality, we get
\begin{align*}
h_k^{\infty}(f)(x) \leq \int_{|x - y| > 1} |f(y)| \left(
\int_0^{+\infty}  \left|t^{k} \frac{\partial^k}{\partial t^k}\, p_t(x, y)
\right|^2 \, \frac{dt}t \right)^{\frac{1}{2}} \, d\mu(y).
\end{align*}

We estimate the inner integral here.

\begin{lem}   \label{P7a}
We have
\begin{align*}   
\left( \int_0^{+\infty}  \left|t^{k} \frac{\partial^k}{\partial t^k}\, p_t(x,
y) \right|^2 \,\frac{dt}t  \right)^{\frac{1}{2}} \lesssim |x-y|^{\frac{k - n}{2}  - \frac{1}{4} }\, e^{-x_1-y_1 -|x-y|}, \qquad |x - y| > 1.
\end{align*}
\end{lem}

\medskip

Since this implies a kernel estimate like  \eqref{glob1} with $q$
replaced by $k + {1}/{2}$,
we can argue as in Section \ref{sec3} using the operator
$\tilde {\mathcal V}_{k+1/2}$. In this way,
 the weak type $(1, 1)$ of $h_1$ and also \eqref{rf2} follow from
Lemma \ref{P7a}.

\medskip

\noindent
{\bf Proof of Lemma \ref{P7a}.} Combining \eqref{6A} and \eqref{KLF2}, we conclude
\begin{align} \label{LLL}
& \int_0^{+\infty}  \left|t^{k} \frac{\partial^k}{\partial t^k}\, p_t(x,
y) \right|^2 \, \frac{dt}t
\lesssim e^{-2 x_1 - 2 y_1 - 2 |x - y|}  \\
&\times \sum |x - y|^{4 p}  \int_0^{+\infty} \left| \frac{|x - y|}{2 t} + 1 \right|^{2 (m - p)}  \left| \frac{|x - y|}{2 t}  - 1 \right|^{2 (m - p)}  t^{2m -4 p - n}\, e^{-2 t (\frac{|x - y|}{2 t}  - 1)^2} \, \frac{dt}t, \notag
\end{align}
where the sum is taken over $m$ and $p$ such that $0 \leq p \leq [\frac{k}{2}]$ and $0 \leq m \leq k - p$. In the last integral, we estimate
\begin{align} \label{72}
\left| \frac{|x - y|}{2 t}  - 1 \right|^{2 (m - p)}  e^{-2 t (\frac{|x - y|}{2 t}  - 1)^2} \lesssim t^{p - m} e^{- t (\frac{|x - y|}{2 t}  - 1)^2}.
\end{align}

Changing variables by $t = s{|x - y|}/{2} $, we see that the integral in the
right-hand side of  \eqref{LLL} is controlled by
\begin{align*}
\int_0^{+\infty} \left( \frac{1}{s} + 1 \right)^{2 (m - p)} \left( |x - y| s \right)^{m - 3 p - n} e^{-\frac12 |x - y| s (\frac{1}{s} - 1)^2} \, \frac{ds}{s}.
\end{align*}
If restricted to the interval  ${1}/{2} < s < 2$, this integral will be of order
of magnitude  $|x - y|^{m - 3 p - n - {1}/{2}}$.
 The other parts of the integral are smaller because of the exponential decay.
From this, Lemma \ref{P7a} follows.    \hfill $\Box$

\medskip

Let us now consider the part at infinity of $g_k,\; k \geq 1$. It is enough to verify the
following estimate:
\begin{align} \label{gk}
\left( \int_0^{+\infty} \left| t^k
\frac{\partial^k}{\partial t^k}\, P_t(x, y) \right|^2 \,\frac{dt}t
\right)^{\frac{1}{2}} \lesssim |x - y|^{-\frac{n+1}{2}}\,
e^{-x_1 - y_1 - |x - y|}, \qquad |x - y| > 1,
\end{align}
since we can then argue as in Section \ref{sec3}.
For this we write the subordination formula as
\begin{align*} 
P_t(x,y) = \int_0^{+\infty}
\phi\left(\frac{t}{2 \sqrt{u}}\right)\,p_u(x,y) \, \frac{du}{u}\,  \qquad \mbox{with} \qquad
\phi(s) = \frac{s}{\sqrt{\pi}}\, e^{- s^2}
\end{align*}
and use Minkowski's integral inequality, getting
\begin{align*}
& \left( \int_0^{+\infty} \left| t^k\, \frac{\partial^k}{\partial t^k}\, P_t(x, y) \right|^2 \, \frac{dt}{t} \right)^{\frac{1}{2}} \\
&\leq \int_0^{+\infty}  \left\{ \left( \int_0^{+\infty} \left| \left( \frac{t}{2 \sqrt{u}} \right)^k \phi^{(k)}\left(\frac{t}{2 \sqrt{u}}\right) \right|^2 \frac{dt}{t} \right)^{\frac{1}{2}} \right\} \, p_u(x, y) \, \frac{du}{u}  \\
&\lesssim \int_0^{+\infty}\, p_u(x, y) \, \frac{du}{u}.
\end{align*}
In view of \eqref{KLF1}, the last expression is majorized by constant times
\begin{align*}
 e^{-x_1 - y_1} \int_0^{+\infty} u^{-\frac{n}{2}}\, e^{- u - \frac{|x -y|^2}{4u}} \, \frac{du}{u}.
\end{align*}
By Lemma \ref{intest} (b), we see that for $|x - y| > 1$ this quantity is no larger than constant times $e^{-x_1 - y_1-|x - y|}\,|x - y|^{-{(n+1)}/{2}}$,
and  \eqref{gk} follows.

\medskip

The next lemma will imply the weak type $(1,1)$ of $H_1$ and also \eqref{rf3}.

\begin{lem} \label{lemma6}
For $|x - y| > 1$, \hskip3pt $t>0$ and $k\ge 1$,
\begin{align*}
\left| t^k \frac{\partial^k}{\partial t^k} \, p_t(x, y) \right| \lesssim
e^{-x_1 - y_1 - |x - y|}\, |x - y|^{\frac{k - n}{2}}.
\end{align*}
\end{lem}

\noindent {\bf Proof.}  From \eqref{6A} and \eqref{KLF2}, we obtain
\begin{align*}
& \left| t^k \frac{\partial^k}{\partial t^k} \, p_t(x, y) \right| \\
&\lesssim e^{- x_1 - y_1 - |x - y|}\, \sum \left| \frac{|x - y|}{2 t} - 1 \right|^{m - p}  \left( \frac{|x - y|}{2 t} + 1 \right)^{m - p}  \, | x - y|^{2 p}\, t^{m - 2 p - \frac{n}{2}} \,
 e^{- t (\frac{|x - y|}{2 t} - 1)^2},
\end{align*}
the sum taken over $0 \leq p \leq [\frac{k}{2}]$, \hskip4pt $0 \leq m \leq k - p$.
We now use  \eqref{72}, where we first take square roots, and write
again  $t = s{|x - y|}/{2} $. The result will be
\begin{align*}
 \left| t^k \frac{\partial^k}{\partial t^k} \, p_t(x, y) \right|
\lesssim e^{- x_1 - y_1 - |x - y|}\, \sum \left( \frac{1}{s}+1 \right)^{m - p} |x - y|^{\frac{p + m - n}{2}} s^{\frac{m - 3 p - n}{2}}\, e^{-\frac{1}{4} |x - y| s (\frac{1}{s} - 1)^2}.
\end{align*}
 Here we simply delete the factor $|x - y|$ in the
last exponent.  Since $p +  m \leq k$, we will get at most the quantity
in the right-hand side of the lemma, multiplied by a bounded function of $s$.

 Lemma \ref{lemma6} is proved.
\hfill $\Box$

\medskip

To verify the sharpness parts of Theorem \ref{TH2}, we use Lemma \ref{lemzz}. For $h_k$ we get with $x \in \Sigma_{\eta}$
\begin{align*}
  \left(h_k\left(\chi_{B(0, 1)}\right)(x)
\right)^2
  \gtrsim \int_{{\eta} ( 1 - 2 {c_1}/{\sqrt{\eta}})/2}^{{\eta} ( 1 -  {c_1}/{\sqrt{\eta}})/2} e^{- 4 \eta}\, \eta^{k - n}\, \frac{dt}{t}
\gtrsim e^{-4 \eta} \,\eta^{k - n - \frac{1}{2}}.
\end{align*}

As in the last part of Section \ref{sec3}, this implies that $h_k$ is not of weak type $(1, 1)$ for $k \geq 2$  and that \eqref{rf2} is sharp.

The case of $H_k$ is even simpler, and this ends the proof of Theorem \ref{TH2}.
\hfill $\Box$

\section*{Acknowledgements}

H.-Q. Li is partially supported by
NSF of China (Grants No. 11571077 and No. 11625102) and ``The Program for
Professor of Special Appointment (Eastern Scholar) at Shanghai
Institutions of Higher Learning''.

\noindent Both authors profited from a grant from the ``Gothenburg Centre for
Advanced Studies in Science and Technology''.


\begin{thebibliography}{00}

\bibitem{A92} G. Alexopoulos, An application of homogenization theory to harmonic analysis: Harnack inequalities and Riesz transforms on Lie groups of polynomial growth, Canad. J. Math. 44 (1992) 691--727.

\bibitem{AS12}  A. Andersson, P. Sj\"ogren, Ornstein-Uhlenbeck theory in finite dimension, Lecture Notes, Dept. Math. Sciences, Chalmers and University of Gothenburg 2012.

\bibitem{An92} J.-P. Anker, Sharp estimates for some functions of the Laplacian on noncompact symmetric spaces, Duke Math. J. 65 (1992) 257--297.

\bibitem{ADY96} J.-P. Anker, E. Damek, C. Yacoub, Spherical analysis on harmonic AN  groups, Ann. Scuola Norm. Sup. Pisa Cl. Sci. (4) 23 (1996) 643--679.

\bibitem{ACDH04} P. Auscher, T. Coulhon, X.-T. Duong, S. Hofmann, Riesz transform on manifolds and heat kernel regularity, Ann. Sci. \'Ecole Norm. Sup. (4) 37 (2004) 911--957.

\bibitem{B87} D. Bakry, \'Etude des transformations de Riesz dans les vari\'et\'es riemanniennes \`a courbure de Ricci minor\'ee, S\'eminaire de Probabilit\'es XXI, Lecture Notes in Math., vol. 1247, Springer, Berlin, 1987, pp. 137--172.

\bibitem{CCH06} G. Carron, T. Coulhon, A. Hassell, Riesz transform and $L^p$-cohomology for manifolds with Euclidean ends, Duke Math. J. 133 (2006) 59--93.

\bibitem{CCFR16} L. Chen, T. Coulhon, J. Feneuil, E. Russ, Riesz Transform for $1 \leq p \leq 2$ Without Gaussian Heat Kernel Bound. J. Geom. Anal. DOI 10.1007/s12220-016-9728-5


\bibitem{CD99} T. Coulhon, X. T. Duong, Riesz transforms for $1 \leq p \leq 2$, Trans. Amer. Math. Soc. 351 (1999) 1151--1169.

\bibitem{CDL03} T. Coulhon, X.-T. Duong, X.-D. Li, Littlewood-Paley-Stein functions on complete Riemannian manifolds for $1 \leq p \leq 2$, Studia Math. 154 (2003) 37--57.

\bibitem{CL04} T. Coulhon, H.-Q. Li, Estimations inf\'erieures du noyau de la chaleur sur les vari\'et\'es coniques et transform\'ee de Riesz, Arch. Math. 83 (2004) 229--242.

\bibitem{ERS99} A.F.M. ter Elst, D. W. Robinson, A. Sikora, Riesz Transforms and Lie Groups of Polynomial Growth, J. Funct. Anal. 162 (1999) 14--51.

\bibitem{FGS94} E. B. Fabes,  C. E. Guti{\'e}rrez, R. Scotto,
Weak-type estimates for the {R}iesz transforms associated with
the {G}aussian measure, {Rev. Mat. Iberoam.} 10  (2)   (1994)    229-281.

\bibitem{GQS92} G. I. Gaudry, T. Qian, P. Sj\"ogren, Singular integrals associated to the Laplacian on the affine group $ax + b$, Ark. Math. 30 (1992) 259--281.

\bibitem{GR07} I. S. Gradshteyn, I. M. Ryzhik, {\em Table of
Integrals, Series, and Products.} 7th edition. Translation edited
and with a preface by Alan Jeffrey and Daniel Zwillinger. Academic
Press, Inc., San Diego, CA, 2007. Reproduction in P.R.China
authorized by Elsevier (Singapore) Pte Ltd.

\bibitem{G09} A. Grigor'yan, {\em Heat kernel and analysis on manifolds.} AMS/IP Studies in Advanced Mathematics, 47. American Mathematical Society, Providence, RI; International Press, Boston, MA, 2009.

\bibitem{HS03} W. Hebisch, T. Steger, Multipliers and singular integrals on exponential growth groups, Math. Z. 245 (2003) 37--61.

\bibitem{Li99} H.-Q, Li, La transformation de Riesz sur les vari\'et\'es coniques, J. Funct. Anal. 168 (1999) 145-238.

\bibitem{LS15} H.-Q. Li, P. Sj\"ogren, Weak type $(1, 1)$ for some  operators related to the Laplacian with drift on real hyperbolic spaces. To appear in Potential Anal. DOI 10.1007/s11118-016-9590-x

\bibitem{LSW14} H.-Q. Li, P. Sj\"ogren, Y.-R. Wu, Weak type $(1, 1)$ of some operators for the Laplacian with drift, Math. Z. 282 (2-3) (2016) 623--633.

\bibitem{L85} N. Lohou\'e, Comparaison des champs de vecteurs et des puissances du laplacien sur une vari\'et\'e riemannienne \`a courbure non positive, J. Funct. Anal. 61 (1985) 164--201.

\bibitem{L87} N. Lohou\'e,  Estimation des fonctions de Littlewood-Paley-Stein sur les vari\'et\'es riemanniennes \`a courbure non positive,  Ann. Sci. \'Ecole Norm. Sup. (4) 20 (1987) 505--544.

\bibitem{L92} N. Lohou\'e, Transform\'ees de Riesz et fonctions sommables, Amer. J. Math. 114 (1992) 875--922.

\bibitem{LM00} N. Lohou\'e, S. Mustapha, Sur les transform\'ees de Riesz sur les espaces homog\`enes des groupes de Lie semi-simples, Bull. Soc. Math. France 128 (2000) 485--495.

\bibitem{LM04} N. Lohou\'e, S. Mustapha, Sur les transform\'ees de Riesz dans le cas du Laplacien avec drift, Trans. Amer. Math. Soc. 356 (2004) 2139--2147.

\bibitem{PS00} S. P{\'e}rez, F. Soria, Operators associated with the {O}rnstein-{U}hlenbeck semigroup, J. London Math. Soc. (2) 61 (3)  (2000) 857--871.

\bibitem{SV08} P. Sj\"ogren, M. Vallarino, Boundedness from $H^1$ to $L^1$ of
Riesz transforms on a Lie group of exponential growth, Ann. Inst. Fourier 58 (2008), 1117--1151.

\bibitem{St70} E. M. Stein, {\em Topics in harmonic analysis related to the Littlewood-Paley theory.}
Ann. of Math. Stud., vol. 63, Princeton Univ. Press, Princeton,
N.J., 1970.

\bibitem{SW69} E. M. Stein, N. J. Weiss, On the convergence of Poisson integrals, Trans. Amer. Math.
Soc. 140 (1969), 35--54.

\end{thebibliography}
\end{document}